
\input amstex
\documentstyle{amsppt}
\baselineskip =8 pt plus 2 pt
\lineskip=8 pt minus 1 pt
\magnification =\magstep1
\lineskiplimit =8pt
\bigskip

\def \pslr{$PSL(2,\Bbb R)$}
\def\G{\Gamma}
\def\g{\gamma}
\def\zo{\overline z}
\def\lg{\Cal L(\Gamma)}
\def\erze{e^r_{\overline{z},\zeta}}
\def\o{\overline}
\def\lo{\text{d}\lambda_0}
\def\lr{\text{d}\lambda_r}
\def\const{\frac{r-1}{\pi}}
\def\ipi{\int_0^{2\pi}}
\def\pslr{PSL(2,\Bbb R)}
\centerline{\bf Some estimates for  Banach space norms
in the von Neumann algebras}
\centerline{\bf associated with the Berezin's
 quantization of  compact Riemann
surfaces}
\bigskip
\centerline{ by Florin R\u adulescu
\footnote {Permanent Address: Department of Mathematics,
 The University of Iowa,
 Iowa City, IA 52246, U.S.A }
\footnote{Member of the Institute of Mathematics,
Romanian Academy, Str. Academiei 14, Bucharest, Romania}}
\bigskip
\proclaim{Abstract}Let $\G$ be any
cocompact, discrete
subgroup  of $\pslr$.
 In this paper we find estimates for the predual and the uniform
Banach space norms  in the von Neumann
 algebras associated with the Berezin' s quantization
of a compact Riemann surface $\Bbb D/\G$. As
a corollary, for large values of the deformation parameter
$1/h$, these von Neumann algebras
are isomorphic.

Using the results in [AS], [AC], [GHJ] on
 the von Neumann dimension of
the Hilbert spaces  in  the  discrete
 series of unitary representations
of $PSL(2,\Bbb R)$, as left modules over $\Gamma$ we deduce that
 the fundamental
 group ([MvN]) of the von Neumann
 $\Cal L(\Gamma)$ contains the positive rational numbers.
Equivalently, this proves
that the algebras $\Cal L(\Gamma)\otimes M_n(\Bbb C)$,
 $n \in \Bbb N$
 are mutually isomorphic.
\endproclaim

In this paper we will find some estimates for
the uniform
(Banach) norm in
  the von Neumann algebras in the Berezin's deformation quantization
of compact Riemann surfaces. We then  use these estimates
together with the results in ([Ra1]) to deduce that the von
Neumann algebras in the Berezin's deformation, are for large
values of the deformation parameter ($r=1/h$), mutually isomorphic.
The results in the  papers ([AS], [Co2], [GHJ]) prove that
the Hilbert spaces, associated with the quantization, are finite,
 left modules over the type $II_1$ factor $\Cal L(\G)$.

Consequently this implies that
 the von Neumann algebras in the deformation
quantization are stably (Morita)
 equivalent with $\Cal L(\G)$. By using the isomorphism result
proved in this paper one deduces
 that the algebras
 $\lg\otimes M_n(\Bbb C)$, $n \in \Bbb N$
 are mutually isomorphic.

This in turn, because of the  existing results
  in the literature, may be used to show
that $\Cal L(\Gamma)$ has (non-irreducible) subfactors having
 non-integer  indices
(by using [Jo1]) or to show that the algebra
 $\Cal L(\Gamma)$ is singly generated
(by using [To]). The above result implies, (by using  [Co1]),
 that there are
type $III_\lambda$ factors whose core is $\lg \otimes B(H)$
for a separable Hilbert space $H$.
 Using [HP], one obtains that there exists no bounded
 projection from $B(H)$ onto the
Banach space subjacent to $\Cal L(\Gamma)$.

 Uniform norm for convolutors in von Neumann algebras of
discrete groups have been first computed by Akemann and Ostrand in
([AO]).
Estimates for the uniform norm in free
group algebras  have been determined by Haagerup in [Ha] and then
used by him to find a non-nuclear $C^\ast$-algebra with the
approximation property.
Some of these estimates have been generalized by Jollisaint ([Jo])
for a larger class of groups (including cocompact groups and
some of Gromov's hyperbolic groups). Some of these results are
used in the papers of Connes and Moscovici ([Co]) on the local
index formulae. We refer to [Pi] (and the references therein) for
the connections between the computation of such norms and
recent developments in the theory of Banach spaces.

The computation for norms of convolutors in
the von Neumann algebra of a free group are part of
the Voiculescu's non-commutative probability theory ([Vo3],
[Vo4]. One very
important consequence of the Voiculescu's theory of random matrices,
as asymptotic models for free group factors, is  that the fundamental
group of the von Neumann algebra of a free group
 with infinitely many generators
contains the rational numbers ([Vo1]).

In this paper
we consider the von Neumann algebras of the Berezin's
quantization of a Riemann surface realized as $\Bbb D/\G$ for
a cocompact subgroup $\Gamma $ of $SU(1,1)$. As we mentioned above,
the algebras in this deformation are stably (Morita) equivalent
 with
$\Cal L(\G)$. In [Ra1] we proved that any element in  the
 von Neumann
algebras in the deformation is represented by a kernel, a
function on $\Bbb D\times\Bbb D$,
$k=k(\o{z},\zeta) $, which is analytic in the second variable
and antianalytic in the first. Moreover the kernel is $\G$
invariant, that is  $k(\o{\g z},\g \zeta)=k(\o{z},\zeta)$,
$\g \in \G$ and $z,\zeta \in \Bbb D$.
We will prove that the uniform norm of the element represented by $k$
is equivalent to the following quantity ($r$ is the reciprocal
of the Planck's constant $h$):
$$\text {max\ }\lbrace\text{sup}_{z \in \Bbb D}
\int_{\Bbb D}|k(\zo,\zeta)| (d(z,\zeta))^r
\text {d}\lambda_0({\zeta}),
\text{sup}_{\zeta\in \Bbb D}
\int_{\Bbb D}|k(\zo,\zeta)| (d(z,\zeta))^r
\text {d}\lambda_0(z)\rbrace.$$

By using this estimate and the results in ([Ra1]) we
prove that for any cocompact,
 fuchsian group subgroup  $\Gamma$
 of the group
$\pslr$ (canonically identified with the group $SU(1,1)$),
the associated von Neumann algebra
 $\Cal L(\Gamma)=\overline{\Bbb C(\Gamma)}^w\subseteq B(l^2(\G))$
has the property that the fundamental group
 $$\Cal F(\Cal L(\Gamma))=
\{t\ |\ \ e(\Cal L(\Gamma)\otimes M_n(\Bbb C))e \cong \Cal L(\Gamma),
 \text{\ for\ some\ projection\ } e
\text {\ of\ trace\ } t \leq n\}$$
 contains $\Bbb Q_{+}\backslash\{0\}$. Equivalently, this proves
that the algebras $\lg\otimes M_n(\Bbb C)$, $n \in \Bbb N$,
 are mutually isomorphic.

As it was pointed out
in [HV] (see also the references therein)
the type  $II_1$ factors associated with cocompact groups in
$PSL(2,\Bbb R)$
have the non-$\Gamma$-property
 of Murray and von Neumann. The property
which we prove in this paper, for the von Neumann algebra
of a cocompact group in $PSL(2,\Bbb R)$,
 is similar to the corresponding property of the
von Neumann algebra of a
free group with infinitely many generators. For this last group,
it was a breakthrough discovery of Voiculescu ([Vo2]), based on
the random matrix model, that the fundamental
 group $\Cal F(\Cal L(F_\infty))$
contains the positive rationals (in fact, as proved in [Ra4], this
model may be used to show that the fundamental group
$\Cal F(\Cal L(F_\infty))$ is $\Bbb R_+\backslash \{0\}$).

The similar
 problem for free groups with finitely many generators is widely
open. Based on Voiculescu's random matrix model for
free groups $F_N$ with finitely many generators, it was
proved independently in ([Dy], [Ra1]) that $\Cal F(\Cal L(F_N))$
is either
$\Bbb R_+\backslash \{0\}$ or either $\{1\}$, independently
on the natural number $N$. The first
situation  would occur if and only if one would have a positive
answer to the von Neumann-Kadison-Sakai question on the
isomorphism of the free group algebras $\Cal L(F_N)$.

We finally note that the only other type $II_1$ factors
(except for the hyperfinite factor) for which
one has some knowledge about the fundamental group, are the
algebras associated with groups with the property $T$. For this
algebras, by a remarkable result of A. Connes, ([Co3]),
we know that the fundamental group is almost countable
(see also [Po] for the recent
 construction of  a different type $II_1$
factor with the same property).

\proclaim{Acknowledgment} During the elaboration of this paper,
I have greatly benefited from disscutions with
C. Frohman, U. Haagerup,
C. T. McMullen, G. Pisier, R. Curto and D. Voiculescu.
\endproclaim

\centerline{\bf Definitions and outline of proofs}
The proof of our main result is based on
some estimates for the Banach space
norms on the von Neumann algebras associated with the
$\Gamma-$ equivariant, Berezin's quantization of
 the unit disk ([Ra2]) or, in other words, for  the algebras
associated with the
quantization deformation ([Be]) of the compact
 Riemann surface  $\Bbb D/\Gamma$.

In this deformation quantization, the  associated
von Neumann
  are type $II_1$ factors ([Ra2], see also
[Ri] for deformations quantization with similar behavior).
 We  prove  that the
 von Neumann algebras corresponding to different values of the
deformation parameter $h$ are all isomorphic, if
 $r=1\backslash h$ is sufficiently large (depending on the
exponent of convergence ([Be], [Pa]) of the group $\Gamma$).
  We refer to
[En] for the computations regarding the asymptotics of the Berezin's
 deformation quantization for such domains
(see also [BMS], [KL], [BC]).

Let $\lambda_r$ be the measure on $\Bbb D$ defined
by $\text{d}\lambda_r(z)=
(1-|z|^2)^{r-2}\text {d} z \text {d} \zo$. Hence $\lambda_0$
is the $PSL(2,\Bbb R)$-invariant measure on $\Bbb D$.
Let $(\pi_r)_{r>1}$ be the (continuous)
  series of projective, unitary
representations of $\pslr$, (identified with $SU(1,1)$),
 on the Hilbert space $H_r=H^2(\Bbb D,\text{d}\lambda_r)$ ([Sa]).
For $M$ a type $II_1$
 factor with normalized trace $\tau$ and for any $t$ in $[0,1]$,
following ([MvN), let
  the reduced algebra $M_t$ be the isomorphism class of the type
$II_1$ factor $eMe$,
 where $e$ is any selfadjoint idempotent of trace $t$ in $M$.

Let   $M_n(\Bbb C)$
carry the canonical (non-normalized) trace.
 The definition for $M_t$
definition also makes sense for any $t>1$ if we replace,
 from the beginning,
the algebra $M$ by $M\otimes M_n(\Bbb C)$, where $n$ is any
 integer bigger then $t$.  The
isomorphism class of $M_t$ is independent of
 the choices made so far ([MvN]).
In particular the fundamental group
 $\Cal F(M)$ is the multiplicative group
$\{t| M_t\cong M\}$.

In [Ra2] we proved that for each $h=1/r>0$
 there exists a suitable vector
space $\Cal V_{h}$ consisting
 of smooth, $\Gamma-$invariant functions on $\Bbb D$ (or simply,
  consisting of smooth functions on $\Bbb D/\Gamma$), so that
$\Cal V_{h}$ is closed under conjugation
 and under the Berezin product
$\ast_h$. Moreover if we endow $\Cal V_h$
 with the trace $\tau$ given by the
integral over a fundamental domain $F$ of
 $\Gamma$ in $\Bbb D$, then,
(by the Gelfand-Naimark -Segal construction), we obtain
 a type $II_1$ factor
$\Cal A_r$ that coincides with the commutant $\{\pi_r(\Gamma)\}'$
 of the image of $\Gamma$ in $B(H_r)$
through the projective, unitary
 representation $\pi_r$ of $\pslr$ into $B(H_r)$.
By $\text{cov\ }\Gamma$ we denote the covolume of $\Gamma$.
 For integer
$r$, (or  if $\Gamma $ is the
 (non-cocompact) group $PSL(2,\Bbb Z)$), it is known,
(see [Co3], [GHJ], [Ra2]), that
the algebras $\Cal A_r$ are isomorphic
to the reduced algebra $\lg_{[(r-1) (\text{cov\ }\Gamma) /\pi]}$.

For $z,\zeta$ in
$\Bbb D$ let $d(z,\zeta)=
(1-|z|^2)^{1/2}(1-|\zeta|^2)^{1/2} |1-\zo\zeta|^{-1}$. This is the
square root of the
hyperbolic cosine of the hyperbolic distance between $z$ and
$\zeta$ in
$\Bbb D$.
For any $z$ in $\Bbb D$ let $e^r_z$ be the
 vector in $H_r$ which corresponds to the
evaluation at $z$. Let  $A$ be a bounded, linear operator
  on $H_r$. Let
$\langle\cdot,\cdot,\rangle_r$ be the scalar product on $H_r$.
 Recall ([Be]) that
the  Berezin's (contravariant) symbol of $A$  is
 a function $\hat A$ on
$\Bbb D\times \Bbb D$, antianalytic in the first variable and
analytic in the second,   computed by
the formula $$\hat A(\zo,\zeta)=\langle A e^r_z,
e^r_\zeta\rangle_r (\langle  e^r_z,e^r_\zeta\rangle_r)^{-1}.$$
 If $A$ commutes with
$\pi_r(\Gamma)$ then the symbol has the following invariance
property: $$\hat A(\gamma \zo,\gamma \zeta)=
\hat A(\zo,\zeta),$$
 for all $z,\zeta$ in $\Bbb D$ and
 for all $\gamma $ in $\Gamma$.

In [Ra2] we introduced the following norm
 ( which is stronger then the uniform norm) on a
weakly dense subalgebra of $B(H_r)$. The
definition of the norm, for $A$ in $B(H_r)$, is  given by
$$||A||_{\lambda,r}=
\text {max\ }\lbrace\text{sup}_{z \in \Bbb D}
\int_{\Bbb D}|\hat A(\zo,\zeta)| (d(z,\zeta))^r
\text {d}\lambda_0({\zeta}),
\text{sup}_{\zeta\in \Bbb D}
\int_{\Bbb D}|\hat A(\zo,\zeta)| (d(z,\zeta))^r
\text {d}\lambda_0(z)\rbrace.$$
Let $\widehat {B(H_r)}$ be the set of all bounded operators on $H_r$
whose $||\cdot||_{\lambda,r}$ norm is finite and let
$\hat{\Cal A}_r$ be
$\Cal A_r\cap\widehat {B(H_r)}$. In the same paper
 ([Ra2]), by using the explicit formulae
for the Berezin's multiplication
 rule $\ast_h$ we determined an explicit formula for
the cyclic, two cocycle $\psi_r$, canonically
 associated with the deformation
(see also [CFS], [CM], [RN]). The formula for $\psi_r$ proved that
$\psi_r$ lives on the algebra $\hat{\Cal A}_r$,
  and that the following estimate
holds
$$ |\psi_r(A,B,C)|\leq
 \text{const}_r ||A||_{\lambda,r}||B||_2 ||C||_2,
\text{\ for \ all\ }A,B, C\in  \hat{\Cal A}_r. \tag 1$$

In general the cyclic cohomology class ([Co]) of the
 cocycle $\psi_r$ represents an obstruction for  the
different products $\ast_h$ to define isomorphic algebras.
The construction of $\psi_r$ may  be used (see [Ra2]) to prove that
the bounded  cohomology  group
  $H^2_{\text{bounded}}(\Gamma, \Bbb Z)$, ([Gr], [Gh]), is
nontrivial for any discrete subgroup
  $\Gamma$ of $\pslr$ having finite covolume.

 If the cocycle $\psi_r$ is bounded,
  by the uniform norm
$||\cdot||_{\infty,r}$ on $\Cal A_r$
  replacing the norm $||\cdot||_{\lambda,r}$
in the equality (1), then standard techniques
([SinS], [CES], [PR]) in the
cohomology theory of
 von Neumann algebras
 are used to show ([Ra2]) that $\psi_r$ is the boundary
of a bounded cycle $\phi_r$.
Hence there exists a bounded, antisymmetric, linear operator
$X_r$ on $L^2(\Cal A_r)$ so that for all
$A,B$ in $L^2(\Cal A_r)$ we have that
 $\phi_r(A,B)=
\langle X_r(A),B\rangle$.
 The evolution operators on $L^2(\Cal A_r)$,
 corresponding to the non-autonomous
differential equation associated with $X_r$ will
  then implement an isomorphism
between the algebras $\Cal A_r$ for different $r$'s ([Ra2]).

In this paper, we show that for a cocompact,
 discrete subgroup $\Gamma $ of $\pslr$, the norms
$||\cdot||_{\lambda,r}$ and $||\cdot||_{\infty,r}$ are equivalent on
$\Cal A_r=\{\pi_r(\Gamma)\}'$. Thus, as  described above,  it
follows that the algebras $\Cal A_r$,
associated with the Berezin's deformation quantization
of $\Bbb D\backslash \G$, are mutually isomorphic.
By ([AS], [Co2], [GHJ])
 the algebras $\Cal A_r=\{\pi_r(\Gamma)\}'$ are
 isomorphic,
for integers $r\geq 2$,
 to $\Cal L(\Gamma)_{(r-1) (\text{cov\ }\Gamma)/\pi})$.

 Note that, (by [Ra2]),
 this  holds in fact
for any $r>1$, if the cocycle
 coming from the projective, unitary representation
$\pi_r|_{\Gamma}$ is trivial in $H^2(\Gamma, \Bbb T)$
 (which happens if e. g.
$\Gamma$ is the (non-cocompact) subgroup $PSL(2,\Bbb Z)$. Hence,
for sufficiently large integers $n,m$, the algebras
$\Cal L(\Gamma)_{[(n-1) (\text{cov\ }\Gamma)/\pi]}$ and
 $\Cal L(\Gamma)_{[(m-1) (\text{cov\ }\Gamma)/\pi]}$ are
isomorphic. This implies that the fundamental group of
$\Cal L(\Gamma)$ contains the positive rational numbers.

Note that if the conjecture in [HV]
  asserting that the von Neumann algebra of a
cocompact, discrete subgroup of $\pslr$ is isomorphic to the  algebra
of  free group whose fractional
  ``number of generators" ([Dy], [Ra1], [Vo1])) depends
 on the covolume of $\Gamma$, then
 it would follow (by [Vo2], [Dy], [Ra1]) that the
question (von -Neumann-Kadison-Sakai,
 [Ka], [Sa]) on the isomorphism of the
algebras $\Cal L(F_N)$ would have an
 affirmative solution. Alternatively
this could happen if one could extend
 the methods in this paper to
non-cocompact groups like $PSL(2,\Bbb Z)$
 or to the discrete Hecke subgroups.

For $z,\zeta$ in
$\Bbb D$ the function
 $d(z,\zeta)= (1-|z|^2)^{1/2}
(1-|\zeta|^2)^{1/2} |1-\zo\zeta|^{-1}$
is the square root of
 the hyperbolic cosines of
 the hyperbolic distance between $z$ and $\zeta$
(see e. g. [Pa]). Denote by
$K_r$ the  symmetric,
 $\Gamma$-equivariant kernel on $\Bbb D$
defined by
$$K_r(z,\eta)=\sum_{\gamma\in \Gamma}
 d(\gamma \eta,z)^r, z,\eta \in \Bbb D.$$
It is well known that the series defining $K_r$
is uniformly convergent
 on compact subsets of $\Bbb D$ if
$r$ is bigger then the double of the exponent of
convergence of the group $\Gamma$ ([Be], [Le], [Pa]). This types
of kernels appear in the Selberg trace formula ([Se]).

 Let
$F$ be any fundamental domain for $\G$ acting on $\Bbb D$.
 Recall that the trace
on $\Cal A_r$ is defined by the formula $\tau_{\Cal A_r}(A)=
\tau(A)= (\lambda_0(F))^{-1} \int_F \hat A(\o{z},z)\lambda_0(z)$
 for any $A$ in $\Cal A_r$
with Berezin symbol $\hat A$. The formula for the product of
two elements $A,B$  in $\Cal A_r$ is
computed out of the symbols $\hat A,\hat B$ as
$$(\hat A\ast \hat B)(\o{z},\zeta)=
c_r\int_{\Bbb D}\frac{\hat A(\o{z},\eta)}{(1-\o{z}\eta)^{-r}}
\frac{\hat B(\o{\eta},\zeta)}
{(1-\o{\eta}\zeta)^r}\text{d}\lambda_0(\eta),z,\zeta\in\Bbb D.$$
Hence the Hilbert space  $L^2(\Cal A_r,\tau)$ associated
 (via the Gelfand-Naimark-Segal) construction
to the trace $\tau$ on the type $II_1$
 factor $\Cal A_r=\{\pi_r(\G)\}'$ is
 identified with the Hilbert space of functions
 $k=k(\o{z},\eta),\g \in\G$ on
$\Bbb D\times\Bbb D$ that are
 antianalytic in the first variable, analytic in the
second which are $\G-$invariant
 ($k(\o{\g z},\g \zeta)=k(\o{z},\zeta)$).
The Hilbert norm is given by the formula:
$$||k||^2_{2,r}=
\int\int_{\Bbb D\times F} |k(\o{z},\eta)|^2
 (d(z,\eta))^{2r}\text{d}\lambda_0(z)
\text{d}\lambda_0(\zeta).$$
Here $\Bbb D\times F$ could be replaced by any fundamental
 domain for the diagonal action of $\G$ on $\Bbb D \times \Bbb D$.

\centerline{\bf The results}

In the next lemma we determine the precise formula for the point
evaluation vectors in the Hilbert space associated to the
deformation quantization. As we pointed out above this may be
identified with a Hilbert space of square summable analytic
 functions and hence it contains evaluation vectors.

\proclaim{Lemma. 1} Let $\Gamma$ be a
cocompact subgroup of $PSL(2,\Bbb R)$.
 For $z,\zeta$ in $\Bbb D$ and for every $r$
  bigger then the double of the exponent of convergence
of $\Gamma$ let $e^r_{\overline{z},\zeta}=
e^r_{\overline{z},\zeta}(\o{\eta_1},\eta_2),
\eta_1,\eta_2 $ in $\Bbb D$ be the function
 on $\Bbb D^2$, antianalytic in the first variable,
analytic in the second defined by the formula:
$$e^r_{\overline{z},\zeta}(\o{\eta_1},\eta_2)=
\frac{r-1}{\pi}\sum_{\gamma}
\frac{(1-\o{\gamma\eta_1}\gamma\eta_2)^r(1-\o{z}\zeta)^r}
{(1-(\o{z})(\g\eta_2))^r(1-\o{\g\eta_1}\zeta)^r}$$
$$=\frac{r-1}{\pi}
\sum_{\gamma}\frac{(1-\o{\eta_1}\eta_2)^r(1-\o{\g z}\g \zeta)^r}
{(1-(\o{\g z})\eta_2)^r(1-(\o{\eta_1})\g\zeta)^r},
 \eta_1,\eta_2\in \Bbb D.$$

Then $\erze$ is the evaluation vector at $z,\zeta\in\Bbb D$
 on $L^2(\Cal A_r,\tau)$, that is
$\tau(A\erze)=\hat A(\o{z},\zeta)$ for all $A$ in $L^2(\Cal A_r,\tau)$.
Moreover $\erze$ belongs to $\hat{\Cal A_r}\subseteq\Cal A_r$.
\endproclaim

Proof. We first prove that
 the series defining $\erze$ is uniformly
convergent on compact subsets in $\Bbb D\times\Bbb D$.
 This will follow automatically
from the computations
 showing that $\erze$ belongs to $\hat{\Cal A_r}$.
We need to estimate
$$\sup_{\eta_1\in\Bbb D}\int_{\Bbb D}
\sum_{\g}\frac {|1-\o{\eta_1}\eta_2|^r|1-\o{\g z}\g \zeta|^r}
{|1-\o{\g z}\eta_2|^r|1-\o{\eta_1}\g \zeta|^r}
 (d(\eta_1,\eta_2))^r\text{d}\lambda_0(\eta_2)$$
$$=\sup_{\eta_1\in\Bbb D}\sum_{\g}
\frac {|1-\o{\g z}\g \zeta|^r}{|1-\o{\eta_1}\g \zeta|^r}
(1-|\eta_1|^2)^{r/2} \int_{\Bbb D}
\frac {1}{|1-\o{\g z}\eta_2|^r}\lambda_{r/2}(\eta_2)$$
$$=\sup_{\eta_1\in\Bbb D}\sum_{\g}
\frac{|1-\o{\g z}\g \zeta|^r(1-|\eta_1|^2)^{r/2}}
{|1-\o{\eta_1}\g \zeta|^r (1-|\g z|^2)^{r/2}}$$
$$=\sup_{\eta_1\in\Bbb D}\sum_{\g}\frac {|1-\o{\g z}\g \zeta|^r}
{(1-|\g z|^2)^{r/2}(1-|\g \zeta|^2)^{r/2}}
\frac{(1-|\eta_1|^2)^{r/2}(1-|\g \zeta|^2)^{r/2}}
{|1-\o{\eta_1}\g \zeta|^r}$$
$$=\sup_{\eta_1\in\Bbb D} (d(z,\zeta))^{-r}\sum_{\g}
\frac{(1-|\eta_1|^2)^{r/2}(1-|\g \zeta|^2)^{r/2}}
{|1-\o{\eta_1}\g \zeta|^r}$$
$$=\sup_{\eta_1\in\Bbb D} (d(z,\zeta))^{-r} K_r(\zeta,\eta_1)
\leq M_r (d(z,\zeta))^{-r}.$$
Hence $\erze$ belongs to $\hat{\Cal A_r}$ and
$$||\erze||_{\lambda,r}\leq M_r (d(z,\zeta))^{-r}.$$

The fact that $\erze$ are the evaluation vectors may
 be tested against elements $A\in \Cal A_r$ which are given as
Toeplitz operators $T^r_{\phi}$ on the Hilbert space
$H^2(\Bbb D,\lambda_r)$ with $\G$-invariant symbol $\phi$.
In this case
$$\tau_{\Cal A_r}(T^r_{\phi}\erze)=
\frac{r-1}{\pi}\int_F\phi(\o{\eta},\eta )
(\sum_{\gamma}\frac{(1-|\gamma\eta|^2)^r(1-\o{z}\zeta)^r}
{(1-\o{z}\g\eta)^r(1-\o{\g\eta}\zeta)^r})\text{d}\lambda_0(\eta)$$
$$=\frac{r-1}{\pi}\int_{\Bbb D}\phi(\o{\eta},\eta )
\frac{(1-\o{z}\zeta)^r}
{(1-\o{z}\g\eta)^r(1-\o{\g\eta}\zeta)^r})\text{d}\lambda_r(\eta)=
<T^r_{\phi}e^r_z,e^r_{\zeta}>.$$

This is exactly the symbol of $T^r_{\phi}$ evaluated at $z,\zeta$.
\proclaim{Remark 2} Estimates for the spectral
distribution of $\erze$, for $z,\zeta\in \Bbb D$
 may  be obtained from
$$||\erze||^2_2=\tau_{\Cal A_r}(\erze e^r_{\o{\zeta},z})=
(\frac{r-1}{\pi})^2
\sum_{\gamma}\frac{(1-\o{\zeta}z)^r(1-\o{\g z}\g \zeta)^r}
{(1-\o{\g z}\zeta)^r(1-\o{\zeta}\g\zeta)^r}, z,\zeta\in \Bbb D.$$
We also note the estimates for the higher moments
 of $\erze e^r_{\o{\zeta},z}$ although
we won't make any use of them. Let $c_r=\frac{r-1}{\pi}$; then
$$\tau_{\Cal A_r}((\erze e^r_{\o{\zeta},z})^n)$$
$$=
(c_r)^{2n}
\sum_{\g_1,...,\g_{2n-1}}
\frac{(1-\o{\g_1 z}\g_1\zeta)^r(1-\g_2z\o{\g_2 \zeta})^r
...(1-\o{\g_{2n-1} z}\g_{2n-1}1\zeta)^r
(1-\o{ z}\zeta)^r}
{(1-z\o{\g_1 z})^r
(1-\g_1 \zeta\o{\g_2 \zeta})^r(1-\g_2 z\o{\g_3 z})^r...
(1-\g_{2n-1}\zeta\o{\zeta})^r}.$$
The following estimate holds  for the
 norm in $L^1(\Cal A_r,\tau)$
of $\erze$:
$$||\erze||_1\leq (\frac{r-1}{\pi})^2
(d(z,\zeta))^{-r}, z,\zeta \in \Bbb D.$$
\endproclaim
Proof. The last estimate may be deduced from the fact that
 the norm of the evaluation vector
$\erze$ in $L^1(\Cal A_r,\tau)$ should be less then
 the norm of the corresponding
(rank 1) evaluation vector (at $z,\zeta$) on $B(H_r)$.
 This norm is easily
computed to be equal to
 $ (\frac{r-1}{\pi})^2(d(z,\zeta))^{-r},$ for all
$ z,\zeta \in \Bbb D.$

\proclaim {Lemma 3} Let $A^\G(\o{\Bbb D}\times \Bbb D)$ be
the space
of (diagonally) $\G$-invariant functions
 on $\Bbb D\times \Bbb D$ that
 are antianalytic in the first variable
and analytic in the second variable, with
the topology of uniform convergence
on compacts subsets of  $\Bbb D\times \Bbb D$.
  For all $\zeta,z$ in $\Bbb D$, the following
integral is convergent in $A^\G(\o{\Bbb D}\times \Bbb D)$ and
is equal to
$e^r_{\o{z},z}$ :
$$\int_{\Bbb D}\erze (d(z,\zeta))^r\text{d}\lambda_0(z)=
e^r_{\o{z},z}.$$
\endproclaim
Proof. We will first check the (absolute)
 uniform convergence on compacts of the integral.
We have for all $\eta_1,\eta_2 \in \Bbb D$ that:
$$\sup_{\eta_1\in\Bbb D}\int_{\Bbb D}
\sum_{\g}\frac {|1-\o{\g\eta_1}\g\eta_2|^r|1-\o{ z} \zeta|^r}
{|1-\o{ z}\g\eta_2|^r|1-\o{\g\eta_1} \zeta|^r}
 (d(z,\zeta))^r\text{d}\lambda_0(\zeta)$$
$$\leq \sum_{\g}\frac {|1-\o{\g\eta_1}\g\eta_2|^r (1-|z|^2)^{r/2}}
{|1-\o{ z}\g\eta_2|^r}\int_{\Bbb D} \frac{(1-|\zeta|^2)^{r/2}}
{|1-\o{\g\eta_1} \zeta|^r}\lo(\zeta)$$
$$\sum_{\g}\frac {|1-\o{\g\eta_1}\g\eta_2|^r (1-|z|^2)^{r/2}}
{|1-\o{ z}\g\eta_2|^r
 (1-|\gamma\eta_1|^2)^{r/2}}$$
$$=\sum_{\g}\frac {|1-\o{\g\eta_1}\g\eta_2|^r}
{(1-|\gamma\eta_1|^2)^{r/2}(1-|\gamma\eta_2|^2)^{r/2}}
\frac{(1-|z|^2)^{r/2}(1-|\gamma\eta_2|^2)^{r/2}}
{|1-\o{ z}(\g\eta_2)|^r}$$
$$=(d(\eta_1,\eta_2))^{-r}K_r(z,\eta_2).$$
This is quantity uniformly bounded for $z,\eta_1,\eta_2$
 in a compact subset of $\Bbb D$ ([Be]).

To check the formula in the statement we need to compute
$$\int_{\Bbb D}\erze(\eta_1,\eta_2)(d(z,\zeta))^r
\text{d}\lambda_0(\zeta)$$
$$\const\int_{\Bbb D}
\sum_{\g}\frac {(1-\o{\g\eta_1}\g\eta_2)^r(1-\o{ z} \zeta)^r}
{(1-\o{ z}\g\eta_2)^r(1-\o{\g\eta_1} \zeta)^r}
 (d(z,\zeta))^r\text{d}\lambda_0(\zeta)$$
$$=\sum_{\g}
\frac {(1-\o{\g\eta_1}\g\eta_2)^r(1-|z|^2)^{r/2}}
{(1-\o{ z}\g\eta_2)^r}(\const)\int_{\Bbb D}
\frac{(1-\o{ z} \zeta)^{r/2}}{(1-\o{\g\eta_1} \zeta)^r}
\frac{1}{(1-z\o{\zeta})^{r/2}}\text{d}\lambda_{r/2}(\zeta)$$
$$=\sum_{\g}\frac {(1-\o{\g\eta_1}\g\eta_2)^r(1-|z|^2)^{r/2}}
{(1-\o{ z}\g\eta_2)^r}\frac{(1-|z|^2)^{r/2}}{(1-\o{\g\eta_1} z)^r}$$
$$=\sum_{\g}\frac {(1-\o{\g\eta_1}\g\eta_2)^r(1-|z|^2)^{r}}
{(1-\o{ z}\g\eta_2)^r(1-\o{\g\eta_1} z)^r}=
e^r_{\o{z},z}(\eta_1,\eta_2).$$
\bigskip
\proclaim{Remark 4} Note that the formula
 in the above statement is related
to the following equality (valid for any $A$ in $\hat {\Cal A_r}$
having the Berezin's symbol $\hat A$). We have
$$(\const)\int_{\Bbb D} A(\o{z},\zeta)(d(z,\zeta))^r
\text{d}\lambda_0(\zeta)=
 A(\o{z},z), z \in \Bbb D.$$
\endproclaim
Proof. Note that the fact that $A$ is in $\hat {\Cal A_r}$
 makes the integral
absolutely convergent. We obtain that
$$(\const)\int_{\Bbb D} A(\o{z},\zeta)(d(z,\zeta))^r
 \text{d}\lambda_0(\zeta)$$
$$=(\const)\int_{\Bbb D}\frac{ A(\o{z},\zeta)(1-|z|^2)^{r/2}}
{(1-\o{ z}\zeta)^{r/2}} \cdot
 \frac{1}{(1-z\o{\zeta})^{r/2}}\text{d}\lambda_{r/2}(\zeta)$$
$$\frac{A(\o{z},z)(1-|z|^2)^{r/2}}{(1-|z|^2)^{r/2}}=A(\o{z},z).$$
This completes the proof.
\bigskip
For a separable Hilbert space $H$ let $\Cal C_1(H)$  denote
the trace class operators on $H$, with the norm $||\cdot||_{1}=
||\cdot||_{1,\Cal C_1(H)}=||\cdot||_1.$
We intend to show that  for $\Gamma$ cocompact subgroup
 of $PSL(2,\Bbb R)$ the
integral in Lemma 3 is also absolutely convergent
 in the normic topology
of $L^1(\Cal A_r,\tau)$. We will first give a formula to estimate
 the norm of an element
in $L^1(\Cal A_r,\tau)$ in terms of its Berezin symbol.
\bigskip

\proclaim{Lemma 5} Let $\G$ be a cocompact subgroup
of $PSL(2,\Bbb R)$. Let $A$ be any element in $L^1(\Cal A_r,\tau)$.
 Let $\g_1,\g_2,...$ be an enumeration of $G$ and let
$F$ be a fundamental domain for $\G$ in $\Bbb D$. Let
$G_N$ be $\cup_{i=1}^N \gamma_i F$.
Let $\chi_{G_N}$ be the characteristic function of
$G_N$ viewed as a multiplication operator on
$L^2(\Bbb D,\lambda_r)$.

 Let
$||\chi_{G_N}A\chi_{G_N}||_{1, \Cal C_1(L^2(G_N,\lambda_r))}=
||\chi_{G_N}A\chi_{G_N}||_{1}$ be the nuclear norm of the
compression of $A$ (viewed as an operator on $L^2(\Bbb D,\lambda_r))$
 to $L^2(G_N,\lambda_r)$.

  Then we have the following formula:
$$||A||_{L^1(\Cal A_r,\tau)}=
\lim_{N\rightarrow \infty}\frac{1}{N}||\chi_{G_N}A\chi_{G_N}||_{1}.$$
\endproclaim
Proof. The normalization $\frac{1}{N}$ comes
 from the fact that the trace of $\chi_{G_N}A\chi_{G_N}$
acting as a (nuclear) operator on $L^2(G_N,\lambda_r)$ is
(by the trace formula in [GHJ]),
$N$ times the trace $$\tau_{\Cal A_r}(A)=
\text{tr}_{B(L^2(F,\lr))}(\chi_F A \chi_F).$$

We use the identification described in [GHJ] (see also
 Chapter 3 in [Ra]) of
$\pi_r(\G)'$ with $\Cal L(\G) \otimes B(L^2(F,\lr))$ acting on
$l^2(\G)\otimes L^2(F,\lr)$. As proved
 in [GHJ] the trace of an element
$x$ in $L^1(\Cal A_r,\tau)$ is computed by the formula
$$\tau_{\Cal A_r}(x)=\text{tr}_{B(L^2(F,\lr))}(\chi_F x\chi_F).$$
In particular
$$||x||_{L^1(\Cal A_r,\tau)}=
\text{tr}_{B(L^2(F,\lr))}(\chi_F |x|\chi_F).$$

 Let
$P_N$ be the projection in $B(L^2(\Bbb D,\lambda_r))$ obtained
by multiplication with the characteristic function of $\chi_{G_N}$.

Denote $M=\{\pi_r(\G)\}'$ and let $x$ be any element in
$L^1(M,\tau)\cap M$ having support a finite projection in $M$
(like the elements in $L^1(\Cal A_r)$ do, as $\Cal A_r=
P_rMP_r$ and $P_r$ is a finite projection in $M$ ([GHJ])).
Then $\chi_{G_N}|x|\chi_{G_N}$ is trace class for any
$N$ and its trace is, (since $|x|$ commutes with $\G$), given by:
 $$N\lbrack(\text{tr}_{B(L^2(F,\lr))}(\chi_F |x|\chi_F)\rbrack.$$

For any element $y$ in $B(L^2(\Bbb D,\lambda_r))$
 denote the positive and negative
 part by $y_+$ and $y_-$ respectively. Then $|P_NxP_N|$,
$(P_N xP_N)_{\pm}$ are weakly convergent to $|x|$, $x_\pm$
Consequently, for any $k=1,2,...$,
$\chi_{\g_k F} (P_N xP_N)_{\pm} \chi_{\g_k F}$
converges weakly to $\chi_{\g_k F} (x_{\pm} )\chi_{\g_k F}$.

The trace of $\frac{1}{N}\chi_{G_N} \lbrack(P_N xP_N)_{\pm}\rbrack
 \chi_{G_N}$
is the same as the trace of the element
 $$A^\pm_N= \frac{1}{N}
 \sum_{k=1}^N\chi_{(\g_k F)} \lbrack(P_N xP_N)_{\pm}\rbrack
 \chi_{(\g_k F)}.$$
The trace of $A^\pm_N$  is in turn equal,
(by bringing back all this elements under
 the projection $P_1=\chi_F$),
to the trace $\text{tr}_{B(L^2(F,\lr))}(B^\pm_N)$
 of the positive element
$$B^\pm_N=\frac{1}{N}\sum_{k=1}^N\chi_F \lbrace\pi_r(\g_k)^\ast
\lbrack \chi_{(\g_k F)} \lbrack(P_N xP_N)_{\pm}
\rbrack \chi_{(\g_k F)}\rbrack (\pi_r(\g_k)
\rbrace\chi_F.$$
On the other hand
since $|x|$ commutes with $\pi_r(\G)$ we have that
$$\frac{1}{N}\sum_{k=1}^N\chi_F \lbrack\pi_r(\g_k)^\ast
( \chi_{(\g_k F)} (x_\pm) \chi_{(\g_k F)}) (\pi_r(\g_k)
\rbrack\chi_F=\chi_F (x_\pm)\chi_F.$$

Since $\chi_{(\g_k F)} ((P_N xP_N)_{\pm}) \chi_{(\g_k F)}$
converges weakly to $\chi_{(\g_k F)}( x_{\pm}) \chi_{(\g_k F)}$
 for any $k$
it follows that $B^\pm_N$ converges
 to $\chi_F (x_\pm) \chi_F.$

Moreover the convergence is dominated;
 all elements are dominated
by  a scalar multiple of the positive
 trace class element $\chi_F P_r\chi_F$ in
$B(L^2(F,\lambda_r))$. Hence,
by Theorem 2.16 in ([Si]), it follows that
$$
\text{tr}_{B(L^2(F,\lr))}(
\frac{1}{N}\chi_{G_N} ((P_N xP_N)_\pm) \chi_{G_N})=
\text{tr}_{B(L^2(F,\lr))}(A^\pm_N)=
\text{tr}_{B(L^2(F,\lr))}(B^\pm_N)$$
converges weakly to $\text{tr}_{B(L^2(F,\lr))}(\chi_F (x_\pm )\chi_F)$.
 This
implies that
$$\lim_{N\rightarrow \infty}
||\frac{1}{N}
\chi_{G_N} x \chi_{G_N})
||_{1,\Cal C_1(L^2(G_N,\lr))}=
\text{tr}_{B(L^2(F,\lr))}(\chi_F |x| \chi_F)$$
$$=
||\chi_F x \chi_F||_{1,\Cal C_1(L^2(F,\lr))}=
||x||_{1,L^1(\Cal A_r)}.$$

To extend  the above result from the class of
 all $x$ in $L^1(M,\tau) \cap M$
having as  support a finite projection in $M$ to the class of  all
 $x$ in $L^1(M,\tau)$ having as
support a finite  projection in $M$ it is sufficient
 to observe the following
inequality.

Note that by the  Pe\-ierls\--Bogo\-liu\-bov inequality,
(also rediscovered by Berezin, see Le\-mma 8.8 in [Si]
 and the references
 therein), we have that, for any positive convex function
$f$ with $f(0)=0$ and $x$ in $M\cap L^1(M,\tau)$,
the following inequality holds true
$$\tau _{\Cal A_r}(f(x))=
\text{tr}_{B(L^2(F,\lr))}(\chi_F f(x) \chi_F)
\geq \text{tr}_{B(L^2(F,\lr))}(f(\chi_F x \chi_F)).$$
Hence for any $N$ and modulo a constant depending on $r$ we have
$$\frac{1}{N}||\chi_{G_N}x\chi_{G_N}||_{1,\Cal C_1(L^2(G_n,\lr))}
\leq \tau _{\Cal A_r}(|x|).$$

This completes the proof.
\bigskip
We mention the following corollary (without proof) since we are
not going to make use of it in this paper. On the other hand
it offers a more tractable (for computations) to obtain estimates
for elements in the predual of $\Cal A_r$.

\proclaim{Corollary} With the notations in Lemma 5 we have that
for any $x$ in $\Cal A_r=\{\pi_r(\G)\}'$ that
$$||x||_{L^1(\Cal A_r,\tau)}
\leq \text{(const)}
\lim \sup_{N\rightarrow \infty}
\frac{1}{\sqrt N}
 ||\chi_{G_N} x\chi_{G_N}||_{2,B(L^2(\Bbb D,\lr))}.$$
\endproclaim

 contains at least
\bigskip
In the next proposition we will use the above estimate to
show that the integral in Lemma 3 is also convergent in
$L^1(\Cal A_r,\tau)$.

\proclaim{Lemma 6}Let $\G$ be a cocompact, discrete subgroup
of $PSL(2,\Bbb R)$. Let $\Cal A_r$ be the von Neumann algebra of
all bounded operators acting on the Hilbert space $H_r$
of the  projective representation $\pi_r$
of $PSL(2,\Bbb R)$,
 that commute with $\pi_r(\Gamma)$.
Then $\Cal A_r$ is a type $II_1$ factor
([GHJ]) and $ L^2(\Cal A_r,\tau)$ is canonically identified
with the Hilbert space of all diagonally $\G-$ invariant functions on
$\Bbb D\times\Bbb D$, antianalytic in the first variable, antianalytic
in the second, which are square summable with respect the
measure $\frac{\lr(z)\lr(\zeta)}{(1-\o{z}\zeta)^{2r}}$, supported on
$\Bbb D\times F$.  For
$z,\zeta \in \Bbb D$  let $\erze \in L^2(\Cal A_r,\tau)$
 be the evaluation
vectors at $z,\zeta$.

 Then, for $r$ sufficiently big,
the integral $$\int_{\Bbb D}
||\erze||_{L^1(\Cal A_r,\tau)}(d(z,\zeta))^r\lo(\zeta)$$ is
absolutely convergent, uniformly in $z$ in
 a compact subset of $\Bbb D$.
\endproclaim
 Before that we insert here a disscution  on
some useful estimates (although they are not of direct use
to the proof itself).

We  evaluate  $||\erze||_{L^1(\Cal A_r,\tau)}$
 for $z,\zeta$ in $\Bbb D$.
Clearly $\erze$ is the sum of
(over $\G$) of the  operators of rank 1 on $L^2(\Bbb D,\nu_r)$ given
by the formula
 $$(1-(\o{\g z})\g\zeta)^r<e^r_{\g z}, \cdot>e^r_{\g\zeta}.$$
It follows that the
nuclear norm
$||\chi_G \erze \chi_G||_{1}$ is bounded by
$$(\const)
 \sum_{\g \in \G} |1-(\o{\g z})\g\zeta|^r
||\chi_G e^r_{\g z}||_{2,L^2(\Bbb D,\lambda_r)}
||\chi_G e^r_{\g\zeta}||_{2,L^2(\Bbb D,\lambda_r)}$$
$$= (\const)\sum_{\g \in \G}|1-(\o{\g z})\g\zeta|^r
\lbrack\int_G\frac{1}{|1-\o{\g z}\eta|^{2r}}\lr(\eta)\rbrack^{1/2}
\lbrack\int_G
\frac{1}{|1-\o{\g \zeta}\eta|^{2r}}\lr(\eta)\rbrack^{1/2}.$$
Hence we have that
$$||\chi_G \erze \chi_G||_{1}$$
$$\leq
\const(d(z,\zeta))^{-r}
\sum_\gamma\lbrack\int_G (d(\g z,\eta)^{2r}\lr(\eta)\rbrack^{1/2}
\lbrack\int_G (d(\g \zeta,\eta)^{2r}\text{d}\lambda_0
(\eta)\rbrack^{1/2}.$$
 Let $\g_1,\g_2,...$ be an enumeration of $\G$ and let
$F$ be a fundamental domain for $\G$ in $\Bbb D$. Let
$G_N$ be $\cup_{i=1}^N \gamma_iF$.
Let $\chi_{G_N}$ be the characteristic function of
$G_N$ viewed as the multiplication operator on
$L^2(\Bbb D,\lambda_r)$. Let $c_r=\const$.
Consequently,
$$\frac{1}{N}\int_{\Bbb D}
 ||\chi_{G_N} \erze \chi_{G_N}||_{1,B(L^2(F,\lambda_r))}
(d(z,\zeta)^r\lo(\zeta)$$
$$\leq
(c_r)\frac{1}{N}\int_{\Bbb D}
\sum_\gamma\lbrack\int_{\cup_{i=1}^N \gamma_i F}
 (d(\g z,\eta)^{2r}\lo(\eta)\rbrack^{1/2}
\lbrack\int_{\cup_{i=1}^N \gamma_i F}
 (d(\g \zeta,\eta)^{2r}\lo(\eta)\rbrack^{1/2}\lo(\zeta)$$
$$=c_r\int_{\Bbb D}
\sum_\gamma\lbrack\frac{1}{N}
\sum_{i=1}^N\int_{\g_i F}(d(\g z,\eta)^{2r}\lo(\eta)\rbrack^{1/2}
\lbrack\frac{1}{N}
\sum_{i=1}^N
\int_{\g_i F}(d(\g\zeta,\eta)^{2r}\lo(\eta)\rbrack^{1/2}\lo(\zeta).$$

We denote the last sums by $\phi_N(z)$.
Assume that $0$ belongs to $F$. It is then easy to see, by the
arguments in ([Le]), that for all $\zeta\in \Bbb D$, $\sigma\in \G$,
 $\int_F(d(\zeta,\sigma\eta)^{2r}\lo(\eta)=
 \int_F(d(\sigma\zeta,\eta)^{2r}\lo(\eta)$ is comparable,
(uniformly
 in $\zeta \in \Bbb D$ and $\sigma \in \Bbb \G$), with $d(z,\sigma 0)$.

Hence (modulo a constant depending on $\Gamma$ and $r$), $\phi_N(z)$
 is dominated
by
$$\sum_\gamma\lbrack\frac{1}{N}\sum_{i=1}^N
 d(\g z, \g_i 0)^{2r}\rbrack^{1/2}
\int_{\Bbb D}
\lbrack\frac{1}{N}\sum_{i=1}^N d(\g \zeta, \g_i 0)^{2r}\rbrack^{1/2}
\lo (\zeta).$$
In turn this quantity
 (because $\lo$ is an invariant measure) is equal to
$$\sum_\gamma
\lbrack\frac{1}{N}\sum_{i=1}^N d(\g z, \g_i 0)^{2r}\rbrack^{1/2}
\int_{\Bbb D}
 \lbrack\frac{1}{N}\sum_{i=1}^N
(d( \zeta, \g_i 0))^{r}\rbrack^{1/2}\lo (\zeta)$$
$$=\int_{\Bbb D}
 \lbrack\frac{1}{N}\sum_{i=1}^N
(d( \zeta, \g_i 0))^{r}\rbrack^{1/2}\lo (\zeta)
\sum_\gamma
\lbrack\frac{1}{N}\sum_{i=1}^N d(\g z, \g_i 0)^{2r}\rbrack^{1/2}.$$
The integral
$\int_{\Bbb D}
 \lbrack\frac{1}{N}\sum_{i=1}^N
(d( \zeta, \g_i 0))^{r}\rbrack^{1/2}\lo (\zeta)$
  is
$$\sum_{\sigma\in\G}
\int_{\sigma F}
 \lbrack
\frac{1}{N}\sum_{i=1}^N d( \zeta, \g_i 0)^{r}
\rbrack^{1/2}\lo (\zeta)$$
$$=\sum_{\sigma\in\G}
\int_{ F} \lbrack\frac{1}{N}\sum_{i=1}^N
d( \sigma\zeta, \g_i 0)^{r}\rbrack^{1/2}\lo (\zeta).$$
Modulo a constant, the last integral
 is comparable (uniformly in $\g_1,\g_2,... $ and $N$)
to
$$\sum_{\sigma\in\G}
 (\frac{1}{N}\sum_{i=1}^N d( \sigma 0, \g_i 0)^{r})^{1/2}.$$
Hence, we get that (modulo a constant depending only on
$\Gamma$ and $r$)
$$
\frac{1}{N}\int_{\Bbb D}
||\chi_{G_N} \erze \chi_{G_N}||_{1,\Cal C_1(L^2(\Bbb D,\lambda_r))}
(d(z.\zeta))^r\lo(\zeta)
$$
$$ \leq
\text{const}_{r,\G}\lbrace\sum_\gamma\lbrack\frac{1}{N}
\sum_{i=1}^N d(\g 0, \g_i 0)^{2r}\rbrack^{1/2}\rbrace^2.$$
Let
$$y_N=\sum_\gamma\lbrack\frac{1}{N}
\sum_{i=1}^N d(\g 0, \g_i 0)^{2r}\rbrack^{1/2},\ N\in\Bbb N.
$$
To be able to take $N$ to limit one should have that the
above sums are uniformly bounded in $N$.
Note that if the summand for $\Gamma$
 wouldn't be raised to the power 1/2 then
 this would have been
(by $\G$ invariance)
$$\sum_\gamma(\frac{1}{N}\sum_{i=1}^N d(\g 0, \g_i 0)^{2r})
=\sum_\gamma(\frac{1}{N}\sum_{i=1}^N d(\g_i\g 0,  0)^{2r})$$
$$=N\sum_\gamma(\frac{1}{N}\sum_{i=1}^N d(\g 0,  0)^{2r})=
\sum_\gamma( d(\g 0,  0))^{2r})$$
which is finite by the arguments in [Be] as
 soon as $k$ is bigger then 1.

We use the notations and the methods in the survey article by Lehner.
Let $n(r,0)$ be the numbers of points in the orbit of $\G 0$
 contained in the
 euclidian disk of radius $r$ in $\Bbb D$. By Tsuji estimates
([Ts]) and by the
asymptotic formula of Huber ([Hu]) $n(r,0)$
 is asymptotically
$\frac{1}{2g-1}\frac{1}{1-r}$,  where $g$ is the genus of
 the compact Riemann surface
$\Bbb D/\G$. Also the distribution of the orbit $\G 0$ is uniform
with respect to arc measure ([EM]). I am very indebted
 to C. T. McMullen
for giving me this information (and many other informations
 which in the end
weren't directly related to this paper).

Neglecting the cardinality of the stabilizer of 0, which
 is finite, and
by using the arguments in [Le]
(in the argument of Theorem 2.2.5 loc. cit) we get

\proclaim{Remark}
Let $\g_1,\g_2,...$ is an enumeration of $\G$ so
 that $d(0,\g_i0)^{-1}$ is
increasing. Let $N=n(s_0,0)$ and let $y_N$ be defined as in
(2) by
$$y_N=\sum_\gamma\lbrack\frac{1}{N}
\sum_{i=1}^N d(\g 0, \g_i 0)^{2r}\rbrack^{1/2},\ N\in\Bbb N.$$
 Then, modulo a constant, we have that
$y_N$ is asymptotically equal to
$$\int_{0}^{1}\int_0^{2\pi}
\lbrack n(s_0,0)^{-1}\int_{0}^{s_0}\int_0^{2\pi}
(\frac {(1-r)^k(1-s)^k}{ (1-rs\exp i(\phi-\theta))^{2k}}
\text{d}\phi \text{d}(n(s,0))\rbrack^{1/2}
\text{d}\theta \text{d}(n(r,0))$$
$$=\int_{0}^{1}\lbrack({1-s_0})\int_{0}^{s_0}
\frac {(1-r)^k(1-s)^k}
{ (1-rs)^{2k-1}}\text{d}(n(s,0))\rbrack^{1/2}\text{d}(n(r,0)).$$

\endproclaim

It is easy to see that, by using the fact that the distribution
$\text{\ d}(n(r,0))$
 is asymptotically $\text{\ const\ }\frac{1}{r-1}$, that
if we add factor $1/\sqrt{N}$ in the sum in (2) or if we add a factor
$\sqrt{(1-s_0)}$ in front of the integral representing the sum then
we would be able to find a finite
 upper bound which is valid for all $N$.

This is because the integral representing the sums is
dominated by terms of the form $(1-s_0)^{-1/2}$ (see the
computations bellow).
 To get rid of this (unfortunate) power of $N$ in our estimate
we have thus to use a better estimate for the integral in our
statement.
\bigskip

Proof (of Proposition 6). We neglect the cardinality of any
stabilizer (because these are finite ([Le])).
We will use the distribution function $n(r,\theta)$ counting the
number of points from the orbit $\Gamma 0$ which are in a sector
of radius $r$ and angle $\theta$ from the origin.
Let
$$\gamma=\pmatrix a&b\\ \o{b}&\o{a}\endpmatrix;\
\gamma_1=\pmatrix a_1&b_1\\ \o{b_1}&\o{a_1}\endpmatrix,$$
$$\gamma_1^{-1}=\pmatrix \o{a_1}&-b_1\\-\o{b_1}&
a_1\endpmatrix\;\ \gamma\gamma_1=
\pmatrix  aa_1+b\o{b_1}&
ab_1+b\o{a_1}\\
\o{b}a_1+\o{a}\o{b_1}&\o{b}b_1+\o{a}\o{a_1}
\endpmatrix.$$
Let $\gamma 0=r\exp(i \theta)=\frac {b}{\o{a}}$,
$\gamma_10=s\exp(i \phi)=\frac {b_1}{\o{a_1}}$. Hence
$$(\gamma\gamma_{1}) 0=
\frac{ab_1+b\o{a_1}}{\o{b}b_1+\o{a}\o{a_1}}=
\frac{\frac {b_1}{\o{a_1}}+ \frac {b}{\o{a}}\frac{\o{a}}{a}}
{\frac {b_1}{\o{a_1}}\frac {\o{b}}{a}
\frac {a}{\o{a}}+1}\cdot
\frac{a\o{a_1}}{\o{a_1}\o{a}}=
\frac
{s\text{e}^{i(\phi+\alpha(\g))}+r\text{e}^{i\theta}}
{1+rs \text{e}^{i(\phi+\alpha(\gamma)-\theta)}}.$$
We will use  the notation:
$$\exp(i\alpha(\gamma))=\exp(i\alpha(\gamma 0))=\frac {a}{\o{a}},
\text{\ if\ } \gamma=\pmatrix a&b\\ \o{b}&\o{a}\endpmatrix \in\G.$$
Note that $$\g^{-1}0=
\pmatrix \o{a}&-b\\-\o{b}&a\endpmatrix 0=
\frac {b}{a}=\frac{b}{\o{a}}\frac{\o{a}}{a}=
(\g 0) \text{e}^{-i\alpha(\g)}=r\text{e}^{i(\theta-\alpha(\g))}.$$
Then
$$(\gamma\gamma_{1}) 0=\frac{r\exp(i (\theta-\alpha(\gamma)))
+s\exp(i \phi)}
{1+rs\exp(i( \phi +\alpha(\gamma)- \theta))}
\exp(i\alpha(\gamma)).$$
Also note that in this case
$$(1-|(\g\g_1) 0|^2)=
d(\g\g_10,0)^2=
d(\g_10,\g^{-1}0)^2$$
$$=\frac {(1-|\g_1 0|^2)(1-|\g^{-1} 0|^2)}
{|1-\o{\g_10}(\g^{-1}0)|^2}=
\frac{(1-r^2)(1-s^2)}
{|1+rs\text{e}^{i(-\phi+\theta-\alpha(\g))}|^2}.$$

In what follows we will use the notations
%
$$\aleph(\g 0)=\aleph(r\text{e}
^{i\theta})=\text{e}^{i(\pi +\alpha(\g))}.$$
%
Hence, with $Z=r\text {e}^{i\theta}$, $\zeta=\g_1 0$, we have that
$$\g\g_10=\frac{Z-\aleph(Z) \zeta}{1-\aleph(Z)\o{Z}\zeta}.$$
Note that if denote $\g\g_10=f_{\g_1}(\g 0)$, then we must have that
for all $\g$ in $\G$ that
$$\sigma(f_{\g_1}(\g 0))=f_{\g_1}(\sigma\g 0),$$
i. e. that
$$
\sigma\lbrack \frac {Z-\aleph (Z)\zeta}
{1-\o{Z}\aleph(Z)\zeta}\rbrack=
\frac{\sigma(Z)-\aleph(\sigma(Z))\zeta}
{1-\o{Z}\aleph(\sigma(Z))\zeta},\tag 1$$
for all $Z$ in the orbit $\G  0$ and all $\sigma$ in $\G$.

For functions $f$ on $\Bbb D$,  we have that
$$\sum _{\gamma \in G} (f(\gamma 0))=\sum _{\g\in \G}
(f(\gamma ^{-1} 0)).$$
By using the density distribution $\text {d} n(r,\theta)$ counting
the points the orbit $\G 0$ in a sector in $\Bbb D$ of radius $r$
and angle $\theta$ (see [Le],[EM]) we get that
$$\int_{\Bbb D} (f(Z\aleph(Z))\lo(Z)-
\int_{\Bbb D} f(Z)\lo(Z)$$
tends to zero when the support of $f$ is close to the boundary
of $\Bbb D$ (modulo terms of lower order
with respect to the distance to the boundary).

 Let  $\g,\g_1\in \G$  and denote
$Z=\g 0,\zeta= \g_1 0$. Then
$$1-\g\g_1 0(\o{\g 0})=
1-\o{Z}\frac{Z-\aleph(Z) \zeta}{1-\aleph(Z)\o{Z}\zeta}=
\frac{(1-|Z|^2)}{1-\aleph(Z)\o{Z}\zeta},\tag 2$$
and
$$|1-\g\g_1 0(\o{\g 0})|=\frac{(1-|Z|^2)}{|1-\aleph(Z)\o{Z}\zeta|}.
\tag 3$$
Also we have for all $\eta_1$ in $\Bbb D$ that
$$(1-(\o{\eta_1})\gamma\zeta)^r=
(1-\o{\eta_1}(\g\g_1 0))^r=
\lbrack 1-\o{\eta_1}
\frac{Z-\aleph(Z)\zeta}
{1-\o{Z}\aleph(Z)\zeta}\rbrack^r $$
$$=\frac
{(1-\o{Z}\aleph(Z)\zeta-\o{\eta_1}(Z-\aleph(Z)\zeta))^r}
{(1-\o{Z}\aleph(Z)\zeta)6r}. \tag 4$$
 With this notations ($Z=\gamma 0$, $\zeta=\g_1 0$) we have that
$$
(1-|\g\g_10|^2)=
\frac{(1-|Z|^2)(1-|\zeta|)^2}{|1-\o{Z}\aleph(Z)\zeta|^2}.\tag 5$$
Recall that for arbitrary $z,\zeta$ in $\Bbb D$  we have that
$$\erze(\o{\eta_1},\eta_2)=
|d(\o{z},\zeta)|^{-r}
\sum_{\g \in\G}
\frac {(1-\o{\g z}(\g \zeta))^r}{|1-\o{\g z}(\g \zeta)|^r}
\frac {(1-\o{\eta_1}\eta_2)^r
(1-|\g z|^2)^{r/2} (1-|\g \zeta|^2)^{r/2}}
{(1-\o{\g z}\eta_2)^r (1-\o{\eta_1}\g \zeta)^r},$$
for all $\eta_i\in \Bbb D.$
Note that the factor $|d(\o{z},\zeta)|^{-r}$ which we get in
front of the formula for $\erze$ will be canceled by the
corresponding factor
$|d(\o{z},\zeta)|^{r}$ in the formula for
$||\erze||_{\lambda,r}$.

The formula which we therefore get for
$\erze$ for $z=0$, $\zeta=\g_1 0$ is
$$e^r_{0,\zeta}(\eta_1,\eta_2)=
(|d(0,\zeta)|^{-r})
\sum_{\g \in \G}
\frac {(1-\o{\g 0}(\g\g_1 0))^r}{|1-\o{\g 0}\g\g_1 0|^r}
\frac {(1-\o{\eta_1}\eta_2)^r
(1-|\g 0|^2)^{r/2}(1-|\g\g_1 0|^2)^{r/2}}
{(1-\o{\g 0}\eta_2)^r(1-\o{\eta_1}(\g\g_1 0))^r},$$
for all $\eta_i\in \Bbb D$.
We now use the method in Lehner ([Le]) to express the sum after
$\gamma$ as an integral. We use the formulae (2), (3), (4) (5)
above.

By using the notation $\zeta=\g_10$ and
 $Z=r\text{e}^{i\theta}$  for the
integration variable ($r,\theta$ are the variables for the
density function $\text{d} n(r,\theta)$ which counts the
number of points in the orbit of $\G 0$ in a sector
of radius $r$ and angle $\theta$ in $\Bbb D$), we get
%
%
%
%
$$e^r_{0,\gamma_1 0}(\eta_1,\eta_2)=
e^r_{0,\zeta}(\eta_1,\eta_2)$$
$$=
(|d(0,\zeta)|^{-r})
\int_{\Bbb D}
\frac {(1-\o{\eta_1}\eta_2)^r (1-|Z|^{r})(1-|\zeta|^2)^{r/2}}
{(1-\o{Z}\eta_1)^r(1-\zeta\aleph(Z)\o{Z}-
 \o{\eta_1}(Z-\aleph(Z)\zeta)^r}
\text{d} n(r,\theta).$$

For our estimates we may replace the measure
$\text{d} n(r,\theta)$ by $r\text{d} n(r,\theta)$ which means
that a sufficiently good approximation for our purposes for
$\erze$ will be
$$ G_{\zeta}(\o{\eta_1},\eta_2)=
(|d(0,\zeta)|^{-r})
\int_{\Bbb D}
\frac {(1-\o{\eta_1}\eta_2)^r (1-|\zeta|^2)^{r/2} }
{(1-\o{Z}\eta_2)^r(1-\zeta\aleph(Z)\o{Z}-
 \o{\eta_1}(Z-\aleph(Z)\zeta)^r}
\text {d}\lambda_{r}(Z).\tag 6$$
The resulting formula, has by the invariance property in (1) has the
property that $$G(\o{\g \eta_1},\g \eta_2)=G(\o {\eta_1},\eta_2).$$
If we further particularize to $\eta_1=\sigma_1 0,
\eta_2=\sigma_2 0$ then we get
$$a_{\sigma_1^{-1}\sigma_2}=G(\o{\sigma_1 0},\sigma_2  0)=
G(0,\sigma_1^{-1} \sigma_2 0)$$
$$=
(|d(0,\zeta)|^{-r})
\int_{\Bbb D}
\frac {(1-|\zeta|^2)^{r/2} }
{(1-\o{Z}(\sigma_1^{-1} \sigma_2 0))^r
(1-\zeta\aleph(Z)\o{Z})^r}
\text {d}\lambda_{r}(Z).\tag 7$$

We  use (7) to find estimates on
$|| e^r_{0,\zeta}||_1$ when $\zeta $ tends to the boundary of
$\Bbb D$.
By Proposition 5 we may estimate the norm $||\cdot||_1$ for an element
whose Berezin ($\G-$ invariant) kernel is $k=k(\o{\eta_1},\eta_2)$
by the norm of the operator on $l^2(\Gamma)$ given by the matrix
$$A=(A_{\sigma_1,\sigma_2});
A_{\sigma_1,\sigma_2}=
|d(\sigma_1 0,\sigma_2 0)|^r k(\o{\sigma_1 0},\sigma_2 0).$$

Let $\delta_{\gamma}$ be the left convolutor by $\gamma$ on
$l^2(\Gamma)$.
Hence to estimate $(|d(0,\zeta)|^{r}) || e^r_{0,\zeta}||_1$
 we could use, by (7),
the norm 1 of the following element in the predual
of $\Cal L(\Gamma)$ :
$$ \sum_{\sigma \in \G} (1-|\sigma 0|^2)^{r/2}
\lbrace\int_{\Bbb D}
\frac {(1-|\zeta|^2)^{r/2} }
{(1-\o{Z}\sigma_1^{-1} \sigma_2 0)^r
(1-\zeta\aleph(Z)\o{Z})^r}
\text {d}\lambda_{r}(Z)\rbrace \delta_\g.$$

We observe that if replace $\aleph(Z)$ by $1$  then
  in the above integral
the power of $\zeta$ in the denominator disappears in the integral
after $Z$. Thus if we replace $\aleph(Z)$ by 1 we get an a convergent
integral for
$\int_{\Bbb D}(|d(0,\zeta)|^{r}) || e^r_{0,\zeta}||_1\lr(\zeta)$.

We have proved above that integrals involving $\aleph(Z)$ behave,
for functions whose support tends to the boundary of $\Bbb D$, like
$\aleph (Z)$ tends to 1. The remainder, by making this approximation
just brings an additional order of zero in the following integrals
which estimate the integral over $\Bbb D$ of the first Sobolev norm of
$(|d(0,\zeta)|^{r})  \erze$. This additional  power of zero
will allow to estimate the norm  $||\cdot||_1$ (as, by ([BS], [St]),
 the
Sobolev $(1+\epsilon)$ norm is dominating the norm $||\cdot||_1$.
The proof of Proposition 7 is complete, once we go through the
following Lemma.


\proclaim{Proposition 7}
Let $\sigma_1,\sigma_2,...$ an enumeration of $\G$ and let
$F$ be a fundamental domain for $\G$ in $\Bbb D$ and let
$G_N=\cup_{i=1}^N \sigma_i F$. Then the
following integrals are (absolutely) convergent, uniformly
in $N\in \Bbb N$ and $z$ in $\Bbb D$.

$$\sup_{N\in\Bbb N}\int_\Bbb D \frac{1}{N}
||\aleph_{G_N} e^r_{z,\zeta}\aleph_{G_N}||_{2,\Cal C_2(L^2(G_N,\lr))}
(d(z,\zeta))^r\lo(\zeta)<\infty.\tag 3$$
Moreover the above expression is bounded, (modulo constants
that only depend on $\G$ and $r$), by the following (finite)
quantity:
$$\sup_{N\in\Bbb N}\sum_{\gamma_1}\lbrack
\frac{1}{N^2}\sum_{i,j=1}^N
\vert
\sum_\gamma \vert d(\gamma0,\sigma_i0)\vert^r
 \vert d(\sigma_j0,\gamma\gamma_10)\vert^r\vert^2\rbrack^{1/2}
<\infty.$$
 Note that the kernel representing
the operator $\erze$ on $L^2(G_N,\lr)$ is\break
$(1-\o{\eta_1}\eta_2)^{-r}\erze(\eta_1,\eta_2)$. Let
$f^{r,N}_{\o{z},\zeta}$
be the partial derivative,  after $\eta_1$, of the kernel
representing $\erze$ on $L^2(G_N,\lr)$. Thus
 $$f^{r,N}_{\o{z},\zeta}(\o{\eta_1},\eta_2)=
\aleph_{G_N}(\eta_1)\aleph_{G_N}(\eta_2) \frac{\text{d}}
{\text {d} \eta_1} (1-\o{\eta_1}\eta_2)^{-r}\erze(\o{\eta_1},\eta_2),
\eta_1,\eta_2\in\Bbb D.$$
Then the following integrals are absolutely convergent,
uniformly in $N$ (and $z\in \Bbb D$):
$$\sup_{N\in\Bbb N}\int_\Bbb D \frac{1}{N}
|| f^{r,N}_{z,\zeta}||_{(L^2(G_N,\lr))^2}
(d(z,\zeta))^r\lo(\zeta)<\infty.$$ Similarly, this is bounded,
(modulo constants that only depend on $r$ and $\G$),
uniformly in $N\in\Bbb N$, by
$$\sup_{N\in \Bbb N}\sum_{\gamma_1}\lbrack
\frac{1}{N^2}\sum_{i,j=1}^N
\vert
\sum_\gamma\frac{1}{( 1-\o{(\sigma_i 0)}\g 0)}
\frac {\vert 1-\o{\g 0}(\g\g_10)\vert}{(1-\o{\g 0}(\g\g_10))}
 (\tilde{d}(\o{\sigma_i0},\gamma 0))^{r}
 (\tilde{d}(\o{\gamma\gamma_1 0},\sigma_j 0))^r\vert^2\rbrack^{1/2}
<\infty.
$$

\endproclaim

Proof.
We have
$$\frac{1}{N^2}
||\aleph_{G_N} e^r_{z,\zeta}\aleph_{G_N}||^2_{2,\Cal C_2(L^2(G_N,\lr))}
$$
$$=
\frac{1}{N^2}
\int_{G_N}\int_{G_N}
\vert \sum_\g\frac {(1-\o{\eta_1}\eta_2)^r
(1-\o{\g z}(\g \zeta))^r}
{(1-\o{\g z}\eta_2)^r (1-\o{\eta_1}(\g \zeta))^r}
\cdot
\frac{1}{( 1-\o{\eta_1}\eta_2)^{r}}
\vert^2
\lr(\eta_1,\eta_2)$$
$$=
\frac{1}{N^2}
\sum_{i,j=1}^N
\int_{\sigma_iF}\int_{\sigma_jF}
\vert \sum_\gamma
\vert d(\g z,\g \zeta)\vert^{-r}
\vert d(\g z, \eta_2)\vert^{r}
\vert d(\eta_1,\g \zeta )\vert^{r}\vert^2\lo(\eta_1,\eta_2).$$
If, as we did above, we replace any integral over the
fundamental domain $F$ with the value of the function to be
integrated at $0$ we get that, modulo a
constant,
$$\frac{1}{N^2}
||\aleph_{G_N} e^r_{z,\zeta}\aleph_{G_N}||^2_{2,\Cal C_2(L^2(G_N,\lr))}
$$
is bounded by the following
$$(d(z,\zeta))^{-2r}
\frac{1}{N^2}
\sum_{i,j=1}^N
\vert \sum_\g d(\g z,\sigma_i 0)^{r}
d(\sigma_j 0,\gamma \zeta)^{r}|^2.
$$
The integral in the statement may, by the same arguments
 as in the
comments after the statement of Proposition 6,
be compared, (uniformly in $z$ in a compact set), by
a discrete sum over $\G$. Hence the integral in our
statement is bounded, (modulo a constant which
depends only on $\Gamma$ and which also
depends (continuously) on $r$), by the following sum:
$$\sum_{\gamma_1}\lbrack
\frac{1}{N^2}\sum_{i,j=1}^N
\vert
\sum_\gamma \vert d(\gamma0,\sigma_i0)\vert^r
 \vert d(\sigma_j0,\gamma\gamma_1 0)\vert^r\vert^2\rbrack^{1/2}.$$
 Here, again, we have replaced the
 integrals over $F$ by the value of the
integrand at $0$.
The integral in which we replace
 $(1-\o{\eta_1}\eta_2)^{-r}\erze(\eta_1,\eta_2)$
by its partial derivative with respect to $\eta_1$:
 $$ \frac{\text{d}}
{\text {d} \eta_1} (1-\o{\eta_1}\eta_2)^{-r}\erze(\o{\eta_1},\eta_2),$$
 will be bounded by a  similar sum, which carries the
additional factor:
$$\frac{1}
{(1-\o{\sigma_i 0}(\gamma 0))}.$$

 If we use the next lemma, this  completes the proof
 of Proposition 7.

\proclaim{Lemma} Let $\G$  be a cocompact subgroup of
$SU(1,1)$. Let $d(z,\zeta)$ be the square root of the
hyperbolic distance between two points $z,\zeta$
 in $\Bbb D$, i.e.
$$d(z,\zeta)=
(1-|z|^2)^{1/2}(1-|\zeta|^2)^{1/2} |1-\zo\zeta|^{-1}.$$
Let $\sigma_1,\sigma_2,...$ be an enumeration of $\Gamma$.
Let $\tilde{d}(\o{z},w)=\frac{(1-|z|^2)^{1/2}(1-|w|^2)^{1/2}}
{(1-\o{z}w)},$ for all $z,w \in \Bbb D$.
The following sums
are absolutely convergent, uniformly with $N\in\Bbb N$:
$$\sup_{N\in \Bbb N}\sum_{\gamma_1}\lbrack
\frac{1}{N^2}\sum_{i,j=1}^N
\vert
\sum_\gamma \vert d(\gamma0,\sigma_i0)\vert^r
 \vert d(\sigma_j0,\gamma\gamma_1 0)\vert^r
\vert^2\rbrack^{1/2}<\infty\tag 8$$
The same holds true (uniformly in $N\in \Bbb N$)
 for the following sums
$$\sup_{N\in \Bbb N}\sum_{\gamma_1}\lbrack
\frac{1}{N^2}\sum_{i,j=1}^N
\vert
\sum_\gamma\frac{1}{( 1-\o{(\sigma_i 0)}\g 0)}
\frac {\vert 1-\o{\g 0}(\g\g_10)\vert}{(1-\o{\g 0}(\g\g_10))}
 (\tilde{d}(\o{\sigma_i0},\gamma 0))^{r}
 (\tilde{d}(\o{\gamma\gamma_1 0},\sigma_j 0))^r\vert^2\rbrack^{1/2}
<\infty.
\tag 9
$$
\bigskip
\endproclaim

We replace the sums $\sum_\gamma$ and $\sum_{\gamma_1}$ in (8)
and (9) by
the integral over $\Bbb D$ with respect to the densities
$\text{d}(n(a_2,\theta_2))$ and respectively
$\text{d}(n(a_1,\theta_1))$. Also we replace the
sums $\sum_{\sigma_1}$, $\sum_{\sigma_2}$ by
$\text{d}(n(t_1,\phi_1))$ and $\text{d}(n(t_2,\phi_2))$
respectively. We also let
$N=\frac{1}{1-s}$. We use the notations from the proof of
Proposition 6.
By using the method in ([Le]) we get that the supreme after
$N$ of the sum (8) in the statement is finite if and only if
the supreme over $s$ of the following
 integrals is finite (for the convenience of the notation
we will replace $r$ by $2r$):
$$
\int_0^1\ipi
\lbrack(1-s)^2
\int_0^s\ipi\int_0^s\ipi
\vert
\int_0^1\ipi
\frac{(1-a_2)^r (1-t_1)^r}
{|1-t_1a_2\exp(i(\phi_1-\theta_2))|^{2r}} \cdot$$
$$\frac{(1-t_2)^r
(1-|\gamma\gamma_1 0|^2)^r}
{ |1-t_2(\gamma\gamma_1 0)|^{2r}  }
\frac{\text {d}\theta_2\text{d} a_2}{(1-a_2)^2}\vert^2
\frac{\text{d}\phi_1 \text{d} t_1}{(1-t_1)^2}
 \frac{\text{d}\phi_2\text{d} t_2}{(1-t_2)^2}\rbrack^{1/2}
\frac{\text {d}\theta_1\text{d} a_1}{(1-a_1)^2}.
$$
$$=
\int_0^1\ipi
\lbrack(1-s)^2
\int_0^s\ipi\int_0^s\ipi
\vert
\int_0^1\ipi
\frac{(1-a_2)^r (1-t_1)^r}
{|1-t_1a_2\exp(i(\phi_1-\theta_2))|^{2r}} \cdot$$
$$\frac{(1-t_2)^r
(1-a_1)^r(1-a_2)^r}
{|1+a_1a_2\text{e}^{i(\alpha(\gamma)-\theta_2+\theta_1)}-
t_2\text{e}^{-i\phi_2}
(a_2 \text{e}^{i\theta_2}+a_1
\text{e}^{i(\alpha(\gamma)+\theta_1)})|^{2r} }\cdot$$
$$
\cdot\frac{\text {d}\theta_2\text{d} a_2}{(1-a_2)^2}\vert^2
\frac{\text{d}\phi_1 \text{d} t_1}{(1-t_1)^2}
 \frac{\text{d}\phi_2\text{d} t_2}{(1-t_2)^2}\rbrack^{1/2}
\frac{\text {d}\theta_1\text{d} a_1}{(1-a_1)^2}.
$$
%
%
%
%
The integral
 in which we replace
 $(1-\o{\eta_1}\eta_2)^{-r}\erze(\o{\eta_1},\eta_2)$ by
its partial derivative
 $$
\frac{\text{d}}
{\text {d}\eta_1}  (1-\o{\eta_1}\eta_2)^{-r}\erze (\o{\eta_1},\eta_2)$$
  with respect $\eta_1$,
carries similar terms. The  term which comes from the
derivation will add a factor, in the summands, of the form
$\frac{1}{( 1-\o{\sigma_i0}({\gamma}0))}$, corresponding to
the differentiation of $(1-\o{\eta_1}\zeta)^{-r}$.
In the above integral this will bring an additional
 factor of the form
 $$\frac{1}{(1-t_1a_2e^{i(\phi_1-\theta_2)})}.$$
Our arguments thus allow to estimate  the sum in (9) by
$$
\int_0^1\ipi
\lbrack(1-s)^2
\int_0^s\ipi\int_0^s\ipi
\vert
\int_0^1\ipi
\frac{(1-a_1)^{r} (1-t_1)^{r} (1-t_2)^{r}}
{(1-t_1a_2e^{i(\phi_1-\theta_2)})^{2r+1}}\cdot$$
$$
\cdot\frac {(1-a_2)^{2r-2}}
{\lbrace1-a_1a_2\text{e}^{i(\alpha(\gamma)-\theta_2+\theta_1)}-
t_2\text{e}^{-i\phi_2}
(a_2 \text{e}^{i\theta_2}+a_1
\text{e}^{i(\alpha(\gamma)+\theta_1)})\rbrace^{2r} }\cdot$$
$$
\cdot \text {d}\theta_2\text{d} a_2\vert^2
\frac{\text{d}\phi_1 \text{d} t_1}{(1-t_1)^2}
\frac{\text{d}\phi_2\text{d} t_2}{(1-t_2)^2}\rbrack^{1/2}
\frac{\text {d}\theta_1\text{d} a_1}{(1-a_1)^2}.
$$
\bigskip
By Stinespring's  estimates ([St]) for  the
trace class norms (applied to the terms of the
form $\aleph_{G_N}\erze\aleph_{G_N}$), the  integrals in  (8) and
(9) above
(the one for $\erze$ and the one corresponding to its
partial derivative), will
bound the integral in the statement of Proposition 7.

This holds because of Theorem 2 in ([St]) which shows
that we may choose the same
constant in the estimates bounding the nuclear norms
for operators on the Hilbert spaces $L^2(G_N,\lr)$ by the
$L^2$ norm on $G_N$  of the first derivative of the kernel
 representing the operator (plus the
Hilbert-Schmidt norm). The Stinespring's theorem applies
here because $\lr$ is a finite
measure on the compact space $\Bbb D$.

 Moreover, the renormalization
we have to perform on $\erze$ (i.e to divide
$\erze(\o{\eta_1},\eta_2)$ by
$(1-\o{\eta_1}\eta_2)^r$)  comes from the fact that
if the Berezin kernel of an operator is $k(\o{\eta_1},\eta_2)$ then
this operator is in fact represented by the kernel
$\frac{k(\o{\eta_1},\eta_2)}{(1-\o{\eta_1}\eta_2)^r}$.

We will show bellow that for $r$ sufficiently big, the integrals are
uniformly bounded in $s\in (0,1)$.
The two integrals bounding (8) and (9) both carry terms of the
form
$$\frac{\text{ \ product\ of \ factors}\
 (1-(f(a_1,a_2,t_1,t_2,s))^{a_f}}
{\text{\ product\ of \ factors}\
(1-(g(a_1,a_2,t_1,t_2,s,\theta_i,\phi_i))^{b_g}},$$
where the functions $f,g$ tend to 1 as the parameters tend to 1.
Moreover the factors on the top of the fraction behave
like products of terms of the form
$$(1-a_i)^{\alpha_i}(1-t_i)^{\beta_i}(1-s)^{\delta}$$
while the terms on the bottom of the fraction behave like that
after taking out the angle measures $\phi_i$,
$\theta_i$, $\alpha(\gamma)$.

We use the following convention to denote an integral of the
above form (or an homegenuous sum of such integrals)
by $[A- B]$ where $A$ is the total degree of the factors
on the top, i.e $A$ is the sum
$$A=\sum_{i=1,2}\alpha_i+\sum_{i=1,2}\beta_i+\delta,$$
and similarly for $B$.

At each of the partial stages in the integration process for
the integrals bounding (8) and (9) we will get
 similar
integrals, with one variable from the set $(a_1,t_1,t_2,a_2)$
(or an angle variable) missing.

The effect of the integration with respect
to $\frac{\text {d} a_i}{(1-a_i)^2}$
 and
$\frac{\text {d} t_i}{(1-t_i)^2}$
 is that they transform integrals (or an homogeneous sum)
 of the
form $[A-B]$ into
an homogeneous sum of integrals of the
type  $[A-1-B]$. The effect of the integration with respect to
$\text {d}\theta_1, \text {d}\theta_2$, $\text {d} \phi_1$ and
$\text {d} \phi_2$ is that they transform
integrals (or an homogeneous sum)
 of the
form $[A-B]$ into
an homogeneous sum of integrals of the
type  $[A+1-B]$.

The effect of the integration of the terms in (8) and (9)
 is explained as
follows. By replacing the measures
$\frac{\text{d}a_i\text{d}\theta_i}{(1-a_i^2)^2}$ and
$\frac{\text{d}t_i\text{d}\phi_i}{(1-t_i^2)^2}$ by,
respectively,
$\frac{a_i\text{d}a_i\text{d}\theta_i}{(1-a_i^2)^2}$ and
$\frac{t_i\text{d}t_i\text{d}\phi_i}{(1-t_i^2)^2}$
we don't change the uniform convergence of the integrals.
The measures are now, for each variable, the measure
$\lo$ on $\Bbb D$. The terms to be integrated may be
 represented at each step
as functions of the hyperbolic distance
for a convenient choice of the variables and we may apply Lemma 1
in [Pat1] (see also ([El])).

For example,
the above mentioned lemma shows that if $|B|<|A|$ then the integral
$$\int_{\Bbb D}\frac{ (1-|\eta|)^{2r}}
{|A-B\o{\eta}|^{2r} | 1-\eta\o{w}|^{2r}}\lo(\eta)=
\int_{\Bbb D}\frac{ (1-|\eta|^2)^{2r}}
{|A|^{2r}|1-(B/A)\o{\eta}|^{2r} | 1-\eta\o{w}|^{2r}}\lo(\eta),$$
 is dominated, (modulo a constant and for some
$\epsilon$ as small as we want), by
$$\frac{1}{(|B|^2-|A|^2)^{r} (1-|w|^2)^{r}}
\lbrack \frac {(|B|^2-|A|^2) (1-|w|^2)}{|A-B\o{w}|^2}
\rbrack^{r-\epsilon}.$$

To evaluate the sums in
(9) (which is majorizing (8))
 we have to go through the following process:
We start with a term of the form $[A-(A+1)]$. Integration
by $\text {d} \theta_2$ gives an homogeneous
sum of terms of the form $[A'-A']$.

 The integration by
$\frac{\text {d} a_2}{(1-a_2)^2}$ will yield an homogeneous
sum of terms of the form  \break$[A''-(A''+1)]$. The square  will give
 homogeneous sum of terms of the  form\break
$[A''-(A''+2)]$. The
recursive  integration by the
$\frac{\text {d} t_i\text {d}\phi_i}{(1-t_i)^2} $  will yield
an homogeneous sum of the type $[A'''-(A'''+2)]$. The square root
will give a similar (eventually an infinite convergent sum)
homogeneous sum of the
type $[A'''-(A'''+1)]$.
 The integral with respect to $\text {d}\theta_1$
and the last integral with respect
to $\frac{\text {d} a_1}{(1-a_1)^2}$ will get us
an homogeneous sum of terms of the form $[A^{(4)}-1-A^{(4)}]$ which
means simply a multiple of $\frac{1}{1-s}$ (to which lower
degree terms in $\frac{1}{1-s}$ are to be
 added, e.g $\frac{1}{(1-s)^\alpha},$
$\alpha<1$).

 The final form of the integral (leaving aside the factor $(1-s)$ and
before performing the last integration by the parameter $a_1$)
is an homogeneous sum of terms of the form
$$ \int_0^1 \frac{(1-s)^A (1-a_1)^B}{(1-sa_1)^{A+B+2}}\text{d} a_1.$$
These (hypergeometric) integrals have the leading term $\frac{1}{1-s}$
(compare with example 8, page 297 in [WW]).
 Finally we have to multiply this
by $(1-s)$ (which comes from under the square root). Thus the
supremumum of the
integrals in
in  (9), after $s$ is finite.

The singularity behavior for the integrals when
the parameters are close to 1, may be explained by the
similarity of this integrals with the Appell's double
hypergeometric functions ([Ex]).

This completes the proof of Lemma 7 and hence the proof of
Proposition 6.
\bigskip
We note that by the same method as above,
 the computation being this time considerably
 easier, allows us to show that
the integral in the Remark after Proposition 6, with an additional
factor $(1-s)^{1/2}$ is convergent. Thus we have the following
statement, which is shows a certain mean convergence for the
sums involved in the determination of Beardon's
exponent of convergence. (In fact for group like the free group
one may check this statement directly by using the natural length
function replacing the logarithm of the hyperbolic distance). This
proposition is not needed for our proof, but it shows why the
first estimate we used for the norm
$||\erze||_1$ fails to give convergence of the integral in
Lemma 6.
\proclaim{Corollary}
 Let $\G$  be a cocompact subgroup of
$SU(1,1)$. Let $\g_1,\g_2,...$ be an
enumeration of $\G$.
 Let $d(z,\zeta)$ be the square root of the
hyperbolic distance between two points $z,\zeta$
 in $\Bbb D$. Then the following sums converge, uniformly in
$N\in \Bbb N$:
$$\sum_\gamma\lbrack\frac{1}{N^2}
\sum_{i=1}^N d(\g 0, \g_i 0)^{2r}\rbrack^{1/2},\ N\in\Bbb N.$$
\endproclaim
\bigskip

We now use the result in Proposition 6 to prove the
 following statement which  estimates the uniform norm
on the von Neumann algebras in the Berezin quantization.

\proclaim {Theorem 8}Let $\Gamma$ be
a cocompact subgroup of $\pslr$. For $r>1$ let
$\pi_r$ be the projective, unitary representation of $\pslr$
 (identified with $SU(1,1)$)
 on the
Hilbert space $H_r=H^2(\Bbb D,\text{d}\lambda_r)$.
Let  $A$ be a bounded operator on $H_r$
 commuting with $\pi_r(\Gamma)$.
Let $|| \cdot||_{\lambda,r}$  be the  norm (initially  defined
 on a weakly dense subalgebra of the commutant
$\Cal A_r=\{\pi_r(\Gamma)\}'$.
 The formula for $|| \cdot||_{\lambda,r}$ is
$$||A||_{\lambda,r}=
\text {max\ }\lbrace\text{sup}_{z \in \Bbb D}
\int_{\Bbb D}|\hat A(\zo,\zeta)| (d(z,\zeta))^r
\text {d}\lambda_0({\zeta}),
\text{sup}_{\zeta\in \Bbb D}
\int_{\Bbb D}|\hat A(\zo,\zeta)| (d(z,\zeta))^r
\text {d}\lambda_0(z)\rbrace.$$
 Then there exists a positive
 constant  $M_r$ and a fixed $r_0>0$ so that
for any $r>r_0$ and
 for all $A$ in $\Cal A_r$ we have that $||A||_{\lambda, r}$
is finite and
$$ ||A||_{\infty,r}\leq ||A||_{\lambda, r}
\leq M_r ||A||_{\infty,r}.$$
 Moreover,
 keeping the symbol $\hat A$ fixed, but
varying $r$ in a bounded interval, the constant $M_r$
remains bounded.
\endproclaim
Proof. We only have to apply Proposition 6
 (and its symmetric version
when the r\^ oles of $z$ and $\zeta$ are switched). The constant
$M_r$ is defined by

$$\max\lbrack\sup_{z\in\Bbb D}(\int_{\Bbb D}
||\erze||_{L^1(\Cal A_r,\tau)}(d(z,\zeta))^r\lr(\zeta)),
\sup_{\zeta\in\Bbb D}(\int_{\Bbb D}
||\erze||_{L^1(\Cal A_r,\tau)}(d(z,\zeta))^r\lr(z))\rbrack.$$

 Then, from the definition
of the norm $|| \cdot||_{\lambda,r}$,
we deduce that for any $A$ in $\Cal A_r$, one has that for all
 $z$ in $\Bbb D$
$$\int_{\Bbb D} |A(\o{z},\zeta)|(d(z,\zeta))^r\lr(\zeta)
= \int_{\Bbb D}
|\tau_{\Cal A_r}(A\erze)|(d(z,\zeta))^r\lr(\zeta)$$
$$\leq ||A||_{\infty,r}\int_{\Bbb D}
||\erze||_{L^1(\Cal A_r,\tau)}(d(z,\zeta))^r\lr(\zeta)).$$
This ends the proof

\proclaim{Theorem 9}  Let $\Gamma$ be a cocompact
discrete subgroup of $\pslr$. For $r>1$ let
$\pi_r$ be the projective unitary representation of $\pslr$
 (identified with $SU(1,1)$)
 on the
Hilbert space $H_r=H^2(\Bbb D,\text{d}\lambda_r)$.
Let $\{\pi_r(\Gamma)\}'$ be
the commutant of  $\pi_r(\Gamma)$ in $B(H_r)$.

Note that the
type $II_1$ factor $\{\pi_r(\Gamma)\}'$ is
 the von Neumann algebras
 associated to  the Berezin's deformation quantization
product $\ast_h$ on functions on $\Bbb D/\G$,
 when $h=1/r$ and the trace is the
 integration over a fundamental domain of
$\Gamma$ in $\Bbb D$.

Then, for any $r$, bigger than a fixed $r_0$,
 the type $II_1$ factors $\{\pi_r(\Gamma)\}'$ are
  mutually isomorphic.

\endproclaim
Proof. Indeed in [Ra2] we proved that,
 (for any fuchsian group $\Gamma$),
the cyclic, two cocycle $\psi_r$ associated with the
 deformation ([CFS],[NT], [Ra2])
 has
the following property
 $$\psi_r(A, B, C)\leq c_r||A||_{\lambda,r}||B||_2 ||C||_2$$
for all $A, B, C$ in $\{\pi_r(\Gamma)\}'$.
 Moreover,  the constants $c_r$ may be
chosen uniformly bounded for $r$ in a bounded interval.
 Also, in [Ra2] we proved
that if one may replace in the above estimate
 the norm  $||\cdot||_{\lambda,r}$ with the
uniform norm $||\cdot||_{\infty,r}$ on $B(H_r)$,
 then the cocycle $\psi_r$ is
a coboundary (on the von Neumann algebra). In this case
 the evolution operator
associated with the operator whose coboundary is $\psi_r$
 (by considering the operator as
 a
quadratic form), will implement (by [Ra2]) an isomorphism between the
algebras $\{\pi_r(\Gamma)\}'$. The preceding
 statement completes thus the proof for
a cocompact group $\Gamma$.

\proclaim {Corollary 10} Let $\Gamma$ be a cocompact, discrete
 subgroup of
$\pslr$. Let $\Cal L(\Gamma)$ be the type $II_1$ factor associated with
$\Gamma$. Then the fundamental group $\Cal F(\Cal L(\Gamma))$ contains the
rational numbers.
 Equivalently the algebras $\Cal L(\Gamma)\otimes M_n(\Bbb C)$
are mutually isomorphic.
\endproclaim

Proof. This follows (by the preceding statement) from the computation in
([AS], [Co3], [GHJ]) that $\{\pi_r(\Gamma)\}'$ is isomorphic to
$\Cal L(\Gamma)_{[(r-1) (\text{cov\ }\Gamma)/\pi]}$
 if $r\geq 2$ is an integer (note that if
$r$ is not an integer, this last
 isomorphism will also hold ([Ra2])
 if the  group cohomology
element in $H^2(\Gamma, \Bbb T)$ associated with the projective, unitary
representation $\pi_r|_{\Gamma}$ vanishes).

Consequently, we obtain  that the
 algebra $\Cal L(\Gamma)_{[(n-1) (\text{cov\ }\Gamma)/\pi]}$
 is isomorphic to
the algebra
$\Cal L(\Gamma)_{[(m-1) (\text{cov\ }\Gamma)/\pi]}$,
 for all sufficiently big integers $n,m$.
 The result then
follows from the fact that $\Cal F(\Cal L(\Gamma))$ is a
multiplicative group.
This completes the proof.

\bigskip

 The following observation is related to the method using in proving
Lemma 3. Although this is not related to the subject of this paper we
mention it here as a consequence of the method used in this paper.
It shows that one may generalize Toeplitz operators
with $\G$-invariant symbol to Toeplitz operators whose symbol
is a finite $\G$- invariant measure on $\Bbb D/\G$. This operators
  are also
bounded  and  commute with $\pi_r(\G)$.

Recall that if $\phi$ is a bounded $\G$-invariant function on
$\Bbb D$, then the corresponding Toeplitz operator $T^r_\phi$
 with symbol
$\phi$ is the compression to $H_r=H^2(\Bbb D,\lambda_r)$ of
the operator of multiplication with $\phi$. Clearly
$T^r_\phi$ commutes with $\pi_r(\G)$, so in the terminology we
used in this paper, $T^r_\phi \in \{\pi_r(\G)\}'$.

 In particular, the next result
 shows that there is no positive constant
$c$ so that  $$c||T^r_\phi||=||T^r_\phi||_{\infty,r}
\geq ||\phi||_{\infty},$$
for all $\G$-invariant bounded measurable functions $\phi$ on $\Bbb D$.
If one drops the condition of $\G$-invariance this was known
to Sarason ([Sar]).

This statement is true
  because if the above inequality would hold for some constant
$c$ then it would follow that any element in $\{\pi_r(\G)\}'$
would be a Toeplitz operator with $\G$-invariant, bounded
 measurable symbol. That these operator do not exhaust all
of $\{\pi_r(\G)\}'$ is the content of the next proposition.

\proclaim{Observation}  Let $\G$ be a cocompact
subgroup of $PSL(2,\Bbb R)$ (identified with
$SU(1,1)$) and let $F$ be a compact
 fundamental domain for the action of
$\G$ on $\Bbb D$.
 Let $\nu$ be a finite measure on $F$. We identify $\nu$ with a
$\G$-invariant measure $\tilde\nu$ on $\Bbb D$ (which is no longer
a finite measure). Consider the quadratic form (eventually unbounded)
$<\cdot,\cdot>_\nu$ defined on $H^2(\Bbb D,\lambda_r)$ by
$$<f,f>_\nu=\int_{\Bbb D}|f|^2\text{d}\tilde\nu(z).$$

Then the quadratic form $<\cdot,\cdot>_\nu$ is bounded and
  defines by
a bounded operator $T^r_\nu$ in $\{\pi_r(\G)\}'$ of uniform norm less
 than a universal
constant (depending on $\G$ and $r$) times the norm $|\nu|(F)$ ([Ru])
of the measure $\nu$:
$$||T^r_\nu||_{\infty,r}\leq \text{const}_{r,\G}|\nu|(F).$$

Note that if $\text{d}\nu=\phi\lo$, for
 a bounded, $\G$-invariant function
$\phi$ than the operator $T^r_\nu$ corresponding
to $<\cdot,\cdot>_\nu$
is $T^r_\phi$.
\endproclaim
Proof. With the notations in Lemma 3,
the (Berezin's) symbol $\hat A$ corresponding to the quadratic form
$<\cdot,\cdot>_\nu$ is computed by the formula
$$\hat A(\o{z},\zeta)=
\int_{\Bbb D} \erze(\eta,\eta)\text{d}\nu(\eta).$$
To check that this is a bounded operator is sufficient (by
([Ra2]) to show that the norm $||A||_{\lambda,r}$ is finite.
Thus we have to estimate
$$\int_{\Bbb D}|\hat A(\o{z},\zeta)|(d(z,\zeta))^{r}\lo(z)$$
$$\leq \int_{\Bbb D}\lbrack\int_F
\sum_{\gamma}\frac{(1-|\gamma\eta|^2)^r(1-\o{z}\zeta)^r}
{(1-\o{z}\g\eta)^r(1-\o{\g\eta}\zeta)^r}
\text{d}|\nu(\eta)|\rbrack
(d(z,\zeta))^{r} \text{d}\lambda_0(\zeta)=$$
$$=\int_F\lbrack\sum_\g \frac{(1-|\g \eta|^2)^r (1-|z|^2)^{r/2}}
{|1-(\g \eta) \o{z}|^r}
\int_{\Bbb D}\frac{1}{|1-(\o{\g \eta})\zeta|^r}\text{d}
\lambda_{r/2}(\zeta)\rbrack \text{d}\nu(\eta)$$
$$\leq \int_F\sum_\g (d(z,\g\eta))^{r}\text{d}\nu(\eta)=
 \int_F K_{r}(z,\eta) \text{d}\nu(\eta).$$
By the estimates in [Le] this quantity is uniformly
bounded in $z$, if $r$ is bigger than twice the exponent of
convergence of $\G$.

\centerline {\bf References}
\item{[AO]} C. Akemann and P. Ostrand, Computing norms
in group $C^ast$-algebras, {\it Amer. J. Math.},
{\bf 98}, (1976), 1015-1047.

 \item{[AS]}            M.F. Atiyah, W. Schmidt,
 A geometric construction
of the discrete series for semisimple Lie groups,
 {\it Invent. Math.}, {\bf 42},
(1977), 1-62.

\item {[Be]} A. F. Beardon,  The exponent of convergence
 of Poincar\' e series,
{\it Proc. London Math. Soc.},{\bf 18}, (1968),461-483.

\item{       [Be]}      F. A. Berezin,  General
 concept of quantization, {\it Comm. Math. Phys.}, {\bf 40}

(1975), 153-174.

\item {[BS]} M. Sh. Birman, M. Z. Solomyak,
Estimates for singular numbers of integral operators,
{\it Vestnik Leningrad University}, {\bf 24},
(1969), 21-28.

\item {[Cob]}  C. A. Berger, L. A. Coburn, Heat
 Flow and Berezin-Toeplitz
estimates, {\it American Journal of Math.}, {\bf 110}, (1994).

\item {[BMS]}
M. Bordemann, E. Meirenken, M. Schlicenmaier,
 Toeplitz quantization of Kaehler manifolds and
$gl(N), N\rightarrow \infty$ limits,
 {\it Mannheimer Manuskripte}, {\bf  147},  1993.

\item{[Co1]}       A. Connes, Non-commutative Differential Geometry,
{\it Publ. Math., Inst. Hautes Etud. Sci.}, {\bf 62}, (1986), 94-144.

    \item{[Co2]} A. Connes, Sur la Theorie
 Non Commutative de l'Integration,
Algebres\newline d'Operateurs, {\it Lecture Notes in Math.},
 {\bf 725}, Springer Verlag.

\item{[Co3]}  A. Connes, Un facteur du type $II_1$ avec le groupe
 fondamentale denombrable,
 {\it J.
Operator Theory}, {\bf 4}, (1980), 151-153.

\item{[CM]}      A. Connes, H. Moscovici, Cyclic Cohomology,
the Novikov Conjecture
 and Hyperbolic Groups, {\it Topology},{\bf  29}, (1990),
345-388.

\item{[CS]} A. Connes, D. Sullivan, Quantized calculus on $S^1$
and quasi-Fuchsian groups.

\item {[CFS]} A. Connes, M. Flato, D. Sternheimer,
 Closed star products and Cyclic Cohomology,
{\it Letters in Math. Physics}, {\bf  24}, (1992), 1-12.

\item{[CES] }      E. Christensen, E. G. Effross, A. M. Sinclair,
 Completely Bounded Multilinear Maps
 and $C^*$-algebraic cohomology,
{ \it Invent. Math.}, {\bf 90}, (1987), 279-296.

\item{[Dyk]}       K. Dykema, Free products of
 hyperfinite von Neumann algebras and
free dimension, {\it Duke Math. J.},{\bf  69},(1993), 97-119.

\item{[El]} Elstrodt, J., Die Resolvente zum Eigenwertproblem
der automorphen Formen in der hyperbolischen Ebene, I,
{\it Math. Ann.}, {\bf 208}, (1974),295-330.

\item{[En]} M. Englis, Asymptotics
 of the Berezin Transform and Quantization
on Planar Domains, {\it Duke Math. Journal}, {\bf 79}, (1995),57-76.

\item{[EM]} Eskin, C. McMullen,
 {\it Duke Math. Journal}, {bf 71}, (1994),

\item{[Ex]}H. Exton, Handbook of hypergeometric integrals,
Ellis Horwood Limited, New York, (1978).

\item{[Gh]}    J. Barge,   E. Ghys, Cocycles d'Euler et de Maslov,
{\it Math. Ann.}, {\bf 294}, (1992), 235-265.

\item{[Gr]}       M. Gromov,
 Volume and bounded  cohomology,
{\it Publ. Math., Inst. Hautes Etud. Sci.},
 {\bf 56}, (1982), pp. 5-100.

\item{[GHJ]}       F. Goodman, P. de la Harpe, V.F.R. Jones, Coxeter
Graphs and Towers of Algebras, {\it Springer Verlag, New York, Berlin,
Heidelberg}, 1989.

\item {[HP]} U. Haagerup, G. Pisier, Bounded linear operators between
$C^{\ast}$-algebras, preprint (1993).

\item {[Ha]} U. Haagerup, An example of
a non nuclear C$^\ast$-algebra which has the matrix
approximation property, {\it Inventiones Math},
 {\bf 50}, (1979), 279-293.

\item {[HV]} P. de la Harpe, D. Voiculescu, A problem on the
type $II_1$ factors of fuchsian groups, Preprint.

\item{[Hu]} H. Huber,  Zur analytischen Theorie
 hyperbolischen Raumformen und
Bewegungsgruppen, {\it Math. Ann.}, {\bf 138}, (1959), 1-26.

\item {[Jo]} P. Jolissaint, Sur la C$^\ast$ alg\` ebre
r\'e duite de certains groupes discrets d'isom\' etries hyperboliques,
{\it C. R. Acad. Sci. Paris}, {\bf 302}, (1986), 657-661.

\item {[Ka]}  R. V. Kadison, Open Problems in Operator Algebras,
Baton Rouge Conference, 1960, mimeographed notes.

\item{[KL]}    S. Klimek, A.Leszniewski,
 {\it Letters in Math. Phys.}, {\bf 24},
(1992), pp. 125-139.

\item{[Le]} J. Lehner, Automorphic forms, in Discrete groups
 and automorphic forms,
editor J. Harvey, 73-119.

\item{[MvN]} F. J. Murray, J. von Neumann, On ring of Operators,IV,
 {\it Annals of Mathematics}, {\bf 44}, (1943), 716-808.

\item{[NT]}    R. Nest,  B. Tsygan, Algebraic index
 theorem for families,
{\it Preprint Series, Kobenhavns Universitet}, {\bf 28}, (1993).

\item{[Pa1]} S. J. Patterson, The
 Exponent of Convergence of Poincar\' e series,
{\it Monatshefte fur Mathematik}, {\bf 82}, (1976), 297-315.

\item {[Pa2]} Spectral Theory and Fuchsian groups,
{\it Math. Proc. Camb. Phil. Soc.}, {\bf 81},
(1977), 59-75.

\item {[Pi]} Quadratic forms in unitary operators, Preprint 1995.

\item {[Po1]} S. Popa, Asymptotic freeness, Preprint 1994.

\item {[PR]} S. Popa, F. R\u adulescu
Derivations of von Neumann Algebras into
 the Compact Ideal Space of a
Semifinite Algebra,
 {\it Duke\ Mathematical\ Journal},{\bf 57}, (1988),
485-518.

\item{[Ra1]}      F. R\u adulescu, Random matrices,
 amalgamated free products and
subfactors in free group factors of noninteger index,
 {\it Inv. Math.} {\bf 115}, (1994), 347-389.

\item{[Ra2]}  F. R\u adulescu,$\G$- invariant form of the
  Berezin quantization
of the upper half plane, preprint, Iowa 1995.

\item{[Ra3]}       F, R\u adulescu,On the von
Neumann Algebra of Toeplitz Operators with
Automorphic Symbol, in Subfactors,
Proceedings of the Taniguchi Symposium on Operator
Algebras, edts. H. Araki,
 Y. Kawahigashi, H. Kosaki, {\it  World Scientific,
 Singapore-New Jersey},1994.

\item {[Ra4]}F. R\u adulescu, The fundamental
 group of the von Neumann algebra of a free
group with infinitely many generators is
 $\Bbb R_+\backslash  \{ 0 \}$,
 {\it Journal\ of\ the\ American\ Mathematical\ Society},
 {\bf 5},
(1992), 517-532.

\item{[Ri]}      M. A. Rieffel, Deformation
 Quantization and Operator
Algebras, {\it  Proc. Symp. Pure Math.}, {\bf 51}, (1990), 411-423.

\item {[Ru]} W. Rudin, Real and Complex Analysis, Addison -Wesley.

\item{[Sak]} Sakai, $C^*$ and $W^*$ Algebras, {\it  Springer Verlag,
 Berlin Heidelberg New York}, 1964.

\item{[Sal]}     P. Sally, Analytic
 Continuation of the Irreducible
Unitary Representations of the Universal
 Covering Group,{\it  Memoirs A. M. S.} , (1968).

\item {[Sar]} D. Sarason, personal communication.

\item{[Se]}      A. Selberg, Harmonic Analysis
 and discontinuous groups in weakly
symmetric Riemannian spaces with applications
 to Dirichlet series.{\it J. Indian
Math. Soc. 20}, (1956),  47-87.

\item{[Si]} B. Simon, Trace class ideals and their applications,
{\it London Math. Lecture Notes Series}, {\bf 35}, 1977.

\item {[SincS]} A. M. Sinclair, R. Smith,
 Hochschild cohomology of von Neumann algebras,
{\it London Math. Soc. Lecture Notes}, {\bf 203}, (1995).

\item {[St]} W. F. Stinespring, A sufficient condition for
an integral operator to have a trace,
{\it Journal Reine Angew. Math.}, {\bf 250}, (1958),
200-207.

\item {[To]} D. Topping, Lectures on von Neumann Algebras,
{\it Van Nonstrand}, London, (1972)

\item{[Ts]} M. Tsuji,Potential Theory in Modern Function Theory,
{\it Maruzen}, Tokyo, 1959.

\item{[Vo1]}       D. Voiculescu, Entropy and Fisher's
 information measure
in free probability, II, {\it Inv. Math.}, {\bf 118}, (1994), 420-443.

\item{[Vo2]}       D. Voiculescu, Circular and semicircular
 systems and free product
factors. In Operator Algebras, Unitary Representations,
 Enveloping algebras and Invariant Theory.
{\it Prog. Math. Boston, Birkhauser},{\bf 92},(1990), 45-60.

\item {[Vo3]} D. Voiculescu, Multiplication of certain non-
commuting random variables, {\it J. Operator Theory},
{\bf 18}, (1987), 223-235.

\item {[WW]} E. T. Whittaker, G. N. Watson, A course in
Modern Analysis, {\it Cambridge University Press}, (1984).

\item{[We]}     A. Weinstein, Lecture at the Bourbaki Seminar,
 Paris, June 1994.

\end